\documentclass[11pt]{amsart}
\usepackage{amsmath,epsf,cite}
\usepackage{amssymb,amsthm,tocvsec2}
\usepackage{graphicx}
\usepackage{epstopdf} 
\usepackage{sidecap}
\usepackage{wrapfig}
\usepackage{pict2e}
\usepackage{float}
\usepackage{graphicx}
\usepackage{setspace}
\usepackage{float}
\usepackage{color}
\usepackage{mathrsfs}
\usepackage{epsfig,epsf,latexsym,subfigure}
\usepackage{fancybox}
\usepackage{fancyhdr}
\usepackage{mdwlist}
\usepackage{paralist}
\usepackage{setspace}
\usepackage{tikz}
\usetikzlibrary{shapes,arrows}
\usepackage{enumitem}
\setlist{nolistsep}
\usepackage{indentfirst}
\usepackage{epstopdf}

\catcode`\@=11
\@addtoreset{equation}{section} 
\date{\today}
\date{ }

\numberwithin{equation}{section} \numberwithin{figure}{section}
\numberwithin{table}{section}

\theoremstyle{plain}

\theoremstyle{remark}

\newcommand{\be}{\begin{equation}}
\newcommand{\ee}{\end{equation}}

\newcommand{\bse}{\begin{subequations}}
\newcommand{\ese}{\end{subequations}}

\def\be{\mathbf{e}}

\newcommand{\ben}{\begin{eqnarray}}
\newcommand{\een}{\end{eqnarray}}
\newcommand{\beq}{\begin{equation}}
\newcommand{\eeq}{\end{equation}}
\newcommand{\bea}{\begin{array}}
\newcommand{\eea}{\end{array}}
\newcommand{\bef}{\begin{figure}[H]}
\newcommand{\eef}{\end{figure}}

\begin{document}
\title[]{ Thermodynamically Consistent  Phase Field Models of Multi-component Compressible Fluid Flows}

\author{Xueping Zhao $\diamond$, Tiezheng Qian $\dagger$ and Qi Wang $\ddagger^{\star}$}
\date{\today}
\thanks{$\diamond$
Department of Mathematics, University of South Carolina, Columbia, SC 29208; Email: xzhao@math.sc.edu. \\
$\dagger$
Department of Mathematics, The Hong Kong University of Science $\&$ Technology, Clear Water Bay, Kowloon, Hong Kong; Email: maqian@ust.hk. \\
$\ddagger^{*}$
Department of Mathematics, University of South Carolina, Columbia, SC 29208; Email: qwang@math.sc.edu}

\begin{abstract}
We present a systematic derivation of thermodynamically consistent hydrodynamic phase field models for compressible viscous fluid mixtures using the generalized Onsager principle along with the one fluid multi-component formulation.  By maintaining momentum conservation while enforcing  mass conservation at different levels, we obtain two   compressible models. When the fluid components in the mixture are incompressible, we show that one compressible model reduces to the quasi-incompressible model via a Lagrange multiplier approach. Several different approaches to arriving at the quasi-incompressible model are discussed.  Finally, we conduct a linear stability analysis on all the binary models derived in the paper and show the differences of the models in near equilibrium dynamics.
\end{abstract}

\maketitle

\section{Introduction}
Fluid mixtures are ubiquitous in nature as well as in industrial applications. In a fluid mixture, when fluid components are compressible, the fluid mixture remains compressible.  While in some fluid mixtures,  when each fluid component is incompressible with a constant specific density,  the fluid mixture may not be incompressible when the densities are not equal.  This fluid mixture was named a quasi-incompressible fluid and its thermodynamically consistent model has been derived and applied to various multi-phase fluid flows \cite{LowengrubRSA1998,Li&WangJAM2014,Gong&Z&Y&W2018_SIAM_JSC,Gong&Zhao&Wang2018_SIAM_JSC2}. The fluid mixture is truly incompressible only when all the fluid components are of the same specific density.  For immiscible fluid mixtures, sharp interface models and phase field models can both be used to describe  fluid motions. While for miscible fluid mixtures, sharp interface models are no longer applicable. So, the phase field model becomes a primary platform to describe the fluid motion in the mixture.

Phase field method has been used successfully to formulate models for fluid mixtures  in many applications like in life sciences \cite{ShaoD1, ShaoD2, WiseJTB2008, ZhaoWang-biofilm-datafit3D, AransonPLOSONE2013} (cell biology \cite{KapustinaPLOSCB2015, Najem&GrantPRE2016, ShaoD1, ZhaoWang-cell-divisionII, AransonJRSI2012, JulicherNP2016}, biofilms \cite{Zhao&WangMB2016, ZhaoWang-biofilm-datafit3D, Zhao&WangBMB2016}, cell adhesion and motility \cite{Camley&Zhao&Li&Levine&RappelPRE2017, AransonSoftMatter2014, Najem&GrantPRE2016,  NonomuraPLOSONE2012, ShaoD1}, cell membrane \cite{Aland&Lowengrub&VoigtJCP2014, GavishN1, WangX1, WitkowskiT1}, tumor growth \cite{WiseJTB2008}), materials science \cite{AlandS1, BordenCMAME2012, Chen&Y1994}, fluid dynamics \cite{LowengrubPRE2009, DuLiRyWa052, LowengrubRSA1998, TorabiS1, BertozziSINUMA2000}, image processing \cite{BertozziIEEE, Li&KimCMA2011}, etc. The most widely studied phase field model for binary fluid mixtures is the one for fluid mixtures of  two incompressible fluids of identical densities \cite{LiSh02, Abels2009}. While modeling binary fluid mixtures using phase field models, one commonly uses a labeling or a phase variable (a volume fraction or a mass fraction) $\phi$ to distinguish between distinct material phases. For instance $\phi = 1$ indicates one fluid phase while $\phi = 0$ denotes the other fluid phase in the binary fluid mixture. For fluid mixtures, the interfacial region is tracked by $0 < \phi < 1$. A transport equation for the phase variable $\phi$ along with the conservation equations of momentum, the continuity equation together with necessary constitutive equations constitute the governing system of  equations for the binary fluid mixture.

In a compressible fluid, the total density $\rho$ is a variable and the mass conservation is given by
\ben\bea{l}
\frac{\partial \rho}{\partial t} + \nabla \cdot ({\rho \bf v}) = 0,
\eea  \label{compress_condi}
\een
where $\bf v$ is the velocity field.
In an incompressible fluid, the mass conservation equation \eqref{compress_condi} reduces to
\ben\bea{l}
\nabla \cdot {\bf v} = 0.
\eea  \label{incompress_condi}
\een
since the density $\rho$ is a constant. In the quasi-incompressible model for fluid mixtures, however, any phase field models adopting continuity condition \eqref{incompress_condi} for such a fluid would be questionable in that it may not meet the  consistency condition with the second law of thermodynamics. In these models, the divergence free condition has to be modified to accommodate quasi-incompressibility. A systematic derivation of phase field models for this type fluid mixture of viscous fluids was given by Lowengrub and Truskinovsky using the mass fraction as the phase variable for binary fluid mixtures \cite{LowengrubRSA1998}  as well as Li and Wang using the volume fraction as the phase variable for multi-component fluid mixtures \cite{Li&WangJAM2014}. The derivations were based on the thermodynamic laws, especially, the second law coupled with the additional constraints imposed by the transport equation of the components consistent with the Onsager linear response theory.

As we know, a hydrodynamic model of single phase incompressible fluids can be derived from the corresponding compressible model by imposing an incompressibility constraint. The resultant is called a constrained theory in continuum mechanics. In nature and industrial applications, there are many material systems  comprising of multi-component compressible as well as incompressible components. For instance, in modeling  tissues, there is the issue of   cell proliferation which makes the volume of the material system and mass grow; in tertiary oil recovery, the mixture of $CO_2$ and n-decane are two important compressible fluid components in the gas-oil mixture; and there are many more material systems in real world applications in this category, where the material components are compressible materials.

In this paper, we  derive thermodynamically consistent  compressible phase field model for multi-component fluid mixtures systematically through a variational approach coupled with the generalized Onsager principle \cite{Yang&Li&Forest&Wang2016}. The generalized Onsager principle consists of the Onsager linear response theory and positive entropy production rule \cite{OnsagerL1,OnsagerL2}. Historically, there have been several theoretical frameworks for one to derive thermodynamical and hydrodynamical models for time
dependent dynamics. The Onsager   principle is the one we adopt in this paper. The Onsager maximum entropy production principle based on the Onsager-Matchlup action potential is another approach to deriving models for Hamiltonian and  dissipative systems \cite{OnsagerL1,DoiM1}. Equivalently, the second law of thermodynamics formulated in the form of Clausius-Duhem relation
 is another classical approach to deriving transient dynamical models \cite{Green-Naghdi1991}. The Hamilton least action principle is a classical one for Hamiltonian or conservative systems. The Hamilton-Rayleigh principle is another incarnation of the Onsager maximum entropy principle \cite{GOUIN1990,Gouin&G1999,Gavrilyuk1998,Gavrilyuk1999,DellIsolaetal2009,Auffray2015}. There are also more elaborate GENERIC and Poisson bracket formalism for non-equilibrium theories \cite{Beri94, Ottinger_book, OttingerPRE1997A, OttingerPRE1997B}. These formulations share the commonality in that the non-equilibrium models have a unified mathematical structure consisting of a reversible (hyperbolic) and irreversible (parabolic, dissipative) component in the evolutionary equations. Some of these equations represent conservation laws for the material system such as mass, momentum and energy conservation while others serve as constitutive equations pertinent to the material properties of the material system that the equations describe. The different methods may differ however in how they handle the boundary conditions as well as if one use the dissipation functional or the mobility (or the friction coefficient) to derive the constitutive equations.

 There are two general approaches to describe multiphasic materials. One uses multi-fluid formulation to describe the density and velocity for each phase explicitly \cite{GOUIN1990,Gouin&G1999,Gavrilyuk1998,Gavrilyuk1999,DellIsolaetal2009,Auffray2015}. Another one uses an average velocity, normally the mass average velocity, together with chemical potentials to describe kinematics for each phase. In the latter approach, the average velocity is a measurable hydrodynamic quantity in fluids. For this reason, we choose this approach to formulate our phase field model for multiphasic fluid flows. Since we consider isothermal fluid systems in this paper, we will use the word "multiphase" and" multi-component" interchangeably.

We formulate the hydrodynamic phase field model for compressible fluid of N-fluid components ($N>1$) using the one fluid multi-component formulation \cite{Beri94}. As it is already demonstrated that hydrodynamic models obeys conservation laws do not necessarily satisfy the second law of thermodynamics if the constitutive equations are not derived in a thermodynamically consistent way. The second law or equivalently the Onsager entropy production requirement is thus an additional condition that a well-posed model should satisfy to ensure its well-posedness mathematically. It does not yield additional governing equations for the model. Instead, it does impose an additional constraint on the model and dictates how entropy is produced during the transient dynamical process when the system approaches the steady state.

In this paper, we first derive two classes hydrodynamic phase field models for compressible fluid mixtures using the Onsager principle. After we obtain the "general" compressible models for multi-component fluid mixtures, we hierarchically impose additional "conservation" and/or "incompressibility" conditions to the material system to arrive at constrained, quasi-incompressible theories to show the hierarchical relationship between the compressible model and the constrained models for multi-component fluid mixtures. Through this systematic approach, we demonstrate how one can derive constrained theories via a Lagrange multiplier approach coupled with the generalized Onsager principle, extending the method applied to single phase materials to multi-component material mixtures in the context of one fluid multi-component framework. In the more general compressible model, we enforce global mass conservation so that the model can be used to describe material systems undergoing mass conversion among different components. We then study near equilibrium dynamics of the general models and their various limits through a linear stability analysis.   Note that we derive the models for viscous fluid components in this paper. However, this approach can be readily extended to complex fluids to account for  viscoelastic effects induced by mesoscopic structures in the complex fluid \cite{Yang&Li&Forest&Wang2016}.

The paper is organized as follows. In \S 2, we formulate two classes of hydrodynamic phase field models for the fluid mixture of compressible fluids and a quasi-incompressible model for the fluid mixture of two incompressible fluids with different mass conservation constraints. In \S 3, we generalize the derivation to compressible fluid mixtures of N-components. The non-dimensionalization of the models is carried out in \S4. In \S 5, we discuss near-equilibrium dynamics of the models using a linear stability analysis. We give the concluding remarks in \S 6.

\section{Hydrodynamic phase field models for binary fluid flows}
We present a systematic derivation of thermodynamically consistent hydrodynamic phase field models for binary compressible fluid flows with respect to various  conditions on mass conservation and incompressibility following the generalized Onsager principle \cite{Yang&Li&Forest&Wang2016}

\subsection{Compressible model with the global mass conservation law}

We first consider a mixture of two compressible viscous fluids with density and velocity pairs $(\rho_1, {\bf v}_1)$ and $(\rho_2, {\bf v}_2)$, respectively. We define the total mass of the fluid mixture as $\rho=\rho_1+\rho_2$ and the mass average velocity as $
{\bf v}=\frac{1}{\rho}(\rho_1{\bf v}_1+\rho_2{\bf v}_2).
$
We allow the mass of fluid components to change via conversion, generation, or annihilation at specified rates.
In this general framework, the mass balance equation for each fluid component is given respectively by
\ben\bea{l}
\frac{\partial \rho_i}{\partial t} + \nabla \cdot(\rho_i {\bf v}_i)  = r_i,  \quad i=1,2,
\eea \label{mass_single}
\een
where $r_i$ is the mass conversion/generation/annihilation rate for the ith component. The corresponding momentum conservation equations are given by
\ben\bea{l}
\frac{\partial (\rho_i {\bf v}_i)}{\partial t}  +  \nabla \cdot (\rho_i {\bf v}_i {\bf v}_i)  = \nabla \cdot \sigma_i  + {\bf F}_{i,e} +  r_i{\bf v}_i, \quad i = 1,2,
\eea \label{momentum_single}
\een
where $\sigma_i$ is the viscous stress of the ith fluid component, ${\bf F}_{i,e} $ the extra force of the ith fluid component including the friction force between different fluid components and some elastic forces, and $r_i{\bf v}_i$ the  force due to mass conversion/generation/annihilation in the ith fluid component.

We rewrite the mass conservation equations using the average velocity as follows
\ben\bea{l}
\frac{\partial \rho_i}{\partial t} + \nabla \cdot(\rho_i {\bf v})  = j_i,  \quad i=1,2,
\eea\label{mass_single}
\een
where $j_i=\nabla \cdot (\rho_i({\bf v}-{\bf v}_i))+r_i$ is the excessive production rate of the ith fluid component.

 If we add the mass balance equations \eqref{mass_single} and linear momentum equations \eqref{momentum_single} of all the components, respectively, we  obtain the total mass balance equation and total linear momentum balance equation as follows
\ben\bea{l}
\frac{\partial \rho}{\partial t}+ \nabla \cdot(\rho{\bf v})  = \sum_{i=1}^2 j_i=  \sum_{i=1}^2r_i,  \\
\\
\frac{\partial (\rho {\bf v})}{\partial t}  +  \nabla \cdot (\rho {\bf v} {\bf v})  = \nabla \cdot \sigma^s+{\bf F}_e,
\eea\een
where ${\bf F}_e = \sum_{i=1}^{i=2} ({\bf F}_{i,e} +  r_i{\bf v}_i)$ and $\sigma^s = \sum_{i=1}^{i=2} (\sigma_i - \rho_i({\bf v}_i - {\bf v})({\bf v}_i - {\bf v}))$ is the  stress tensor. The angular momentum balance implies the symmetry of $\sigma^s$. All $j_i$, i=1, 2, $\sigma^s$ and ${\bf F}_e$ will be determined later through constitutive relations.

We assume the free energy of the system is given by
\ben\bea{l}
F = \int_V f(\rho_1, \rho_2, \nabla \rho_1, \nabla \rho_2) d{\bf x},
\eea\een
where f the free energy density function and $V$ the domain in which the fluid mixture occupies.
The total mechanical energy of the system is given by
\ben\bea{l}
E_{total} = \int_V [\frac{1}{2}\rho ||{\bf v}||^2 + f ]d{\bf x} .
\eea\een

We next calculate the total energy dissipation rate as follows.
\ben\bea{l}
\frac{dE_{total}}{dt}
=  \int_{V}[- \sigma^s : {\bf D} +  ({\bf F}_e + \rho_1 \nabla \mu_1 + \rho_2 \nabla \mu_2 - \frac{1}{2}(j_1 + j_2){\bf v}) \cdot {\bf v}   +   \mu_1  (j_1) + \mu_2 (j_2) ] d{\bf x}  \\
+ \int_{\partial V} [(\sigma^s \cdot {\bf v})    - \frac{1}{2} (\rho {\bf v} \|{\bf v}\|^2)
+  ( -  \mu_1 \rho_1 {\bf v} -  \mu_2 \rho_2 {\bf v}  +\frac{\partial f}{\partial (\nabla \rho_1)}\frac{\partial \rho_1}{\partial t}    +   \frac{\partial f}{\partial (\nabla \rho_2)}\frac{\partial \rho_2}{\partial t} ) ]\cdot {\bf n} dS.
\eea
\een
where  $\mu_1 = \frac{\partial f}{\partial \rho_1} - \nabla \cdot \frac{\partial f}{\partial \nabla \rho_1}$, $\mu_2= \frac{\partial f}{\partial \rho_2} - \nabla \cdot \frac{\partial f}{\partial \nabla \rho_2}$ are the chemical potentials with respect to $\rho_1$ and $\rho_2$, respectively,  ${\bf D} = \frac{1}{2}(\nabla {\bf v} + \nabla {\bf v}^T)$ is the rate of strain tensor.
We define the elastic force as
\ben\bea{l}
{\bf F}_e = -\rho_1 \nabla \mu_1 - \rho_2  \nabla \mu_2 + \frac{1}{2}(j_1 + j_2){\bf v}.
\eea\een
This force does not contribute to the energy dissipation.

Using the Onsager principle, we propose
\ben\bea{l}
\sigma^s =  2 \eta {\bf D} +  {\nu} tr({\bf D}){\bf I},\\
\left (
\bea{l}
j_1\\
j_2
\eea \right ) =-{\mathcal M} \left(
\bea{l}
\mu_1\\
\mu_2
\eea \right ),
\eea\label{general}
\een
where $\eta, {\nu}$ are mass-average shear and volumetric viscosities, respectively, and $\mathcal M$ is an operator. The bulk energy dissipation rate reduces to
 \ben\bea{l}
\frac{dE_{total}}{dt}
=  -\int_{V}[ 2\eta {\bf D} :{\bf D} +{\nu} tr({\bf D})^2 +  (\mu_1,\mu_2) \cdot {\mathcal M} \cdot ( \mu_1
,  \mu_2)]  d{\bf x}.
\eea\een
It is non-positive definite provided $\mathcal M$ is nonnegative definite and $\eta, {\nu}$ are non-negative.
The constitutive relation gives a general compressible model for binary fluid flows.

In practice, the interesting scenarios are the following two:
\begin{enumerate}
\item  $\int_V \sum_{i=1}^2 r_i=0$; so, $\int_V \sum_{i=1}^2 j_i=0.$
\vskip 12 pt
\item $r_i=0, i=1,2$; so, $\sum_{i=1}^2 {j}_i=0$.
 \end{enumerate}
 The first condition yields the compressible model of global mass conservation law while the second one gives the compressible model of local mass conservation law. For the first case, a special choice of the mobility operator is the following
\ben\bea{l}
j_1 =  \nabla \cdot M_{11} \nabla \mu_1+\nabla \cdot M_{12} \nabla \mu_2,\\\\
j_2 =  \nabla \cdot M_{21} \nabla \mu_1+\nabla \cdot M_{22} \nabla \mu_2,
\eea\een
where  $M_{ij}, i,j=1,2$ are mobility coefficients.
If we set
\ben\bea{l}
{\bf v}|_{\partial V}=0,  \qquad {\bf n} \cdot \frac{\partial f}{ \partial (\nabla \rho_1)}|_{\partial V} = 0,\qquad {\bf n} \cdot \frac{\partial f}{ \partial (\nabla \rho_2)}|_{\partial V} = 0,
\eea \label{pdes-bc}
\een
on the boundary of the domain $V$, the surface terms vanish in the energy dissipation functional so that the energy dissipation rate reduces to
\ben\bea{l}
\frac{dE_{total}}{dt}
=  -\int_{V}[ 2\eta {\bf D} :{\bf D} +{\nu} tr({\bf D})^2 +  (\nabla \mu_1,\nabla \mu_2) \cdot {\bf M} \cdot (\nabla \mu_1
, \nabla \mu_2)]  d{\bf x},
\eea\een
where ${\bf M}=(M_{ij}).$
It is non-positive definite provided $\eta, {\nu} \geq 0$ and ${\bf M}$ is non-negative definite.

We summarize the governing system of equations in  the  hydrodynamic model for binary compressible fluids with a global mass conservation law as follows:
\ben\bea{l}
\begin{cases}
\frac{\partial \rho_1}{\partial t} + \nabla \cdot (\rho_1 {\bf v}) = j_1 =   \nabla \cdot M_{11} \cdot \nabla \mu_1+\nabla \cdot M_{12} \cdot \nabla \mu_2,\\\\
\frac{\partial \rho_2}{\partial t} + \nabla \cdot (\rho_2 {\bf v}) = j_2 =   \nabla \cdot M_{21} \cdot \nabla \mu_1+\nabla \cdot M_{22} \cdot \nabla \mu_2,\\\\
\frac{\partial (\rho {\bf v})}{\partial t}  +  \nabla \cdot (\rho {\bf v} {\bf v}) - \frac{1}{2}(j_1+j_2){\bf v}  = 2 \nabla \cdot ( \eta{\bf D}) + \nabla ( {\nu} \nabla \cdot {\bf v})-  \rho_1 \nabla \mu_1 - \rho_2 \nabla \mu_2.
\end{cases}
\eea \label{compressible_sun}
\een
We denote the shear viscosities of the fluid component 1 and 2 as $\eta_1,\eta_2$, and the volumetric viscosities of the two components as $\nu_1, \nu_2$, respectively. There are several options of defining average viscosity coefficients in the binary model.
\begin{itemize}
\item  Viscosity coefficients are interpolated using mass fractions and given by
\ben\bea{l}
\eta = \frac{\rho_1}{\rho} \eta_1 + \frac{\rho_2}{\rho} \eta_2,\qquad {\nu} = \frac{\rho_1}{\rho} {\nu}_1 + \frac{\rho_2}{\rho} {\nu}_2.
\eea\een
\item  Viscosity coefficients are interpolated through volume fractions $\phi$ and $(1-\phi)$ in quasi-incompressible models (presented later) and given by
\ben\bea{l}
\eta = \phi \eta_1 + (1-\phi) \eta_2,\qquad {\nu} =\phi {\nu}_1 + (1-\phi) {\nu}_2,
\eea\een
where $\phi$ is the volume fraction of fluid 1.
  \item By the Krieger-Dougherty law, the shear viscosity $\eta$ exhibits a strong non-linear dependence on the local solute concentration and is given by
\ben
\eta(x) = \eta_0 (1-x)^{-\nu},
\een
in which $x$ is the solute concentration ($\rho_1$ or $\rho_2$ in this model), $\eta_0$ is the viscosity of the pure solvent. For example in mixtures of $CO_2$ and n-decane, the solvent is n-decane and solute is $CO_2$. The volumetric viscosity is obtained analogously.
\end{itemize}

As a customary approximation, we assume the free energy density function $f$ is composed of the conformational entropy, and the bulk energy $h$ as follows
\ben\bea{l}
f(\rho_1, \rho_2, \nabla \rho_1, \nabla \rho_2) = h(\rho_1, \rho_2, T) \\+ \frac{1}{2}(\kappa_{\rho_1 \rho_1} (\nabla \rho_1)^2 + 2\kappa_{\rho_1 \rho_2} (\nabla \rho_1, \nabla \rho_2) + \kappa_{\rho_2 \rho_2} (\nabla \rho_2)^2).
\eea \label{eq:kappa_coef}
\een
where $T$ is the absolute temperature, $h(\rho_1, \rho_2, T)$ is the homogeneous bulk free energy density function, $\kappa_{\rho_1 \rho_1}, \kappa_{\rho_1 \rho_2}$ and $\kappa_{\rho_2 \rho_2}$ are parameters parameterizing the conformational entropy, which are all functions of $T$. For example, for the partially miscible binary fluid mixture of n-decane and $CO_2$, where n-decane is denoted as fluid 1 and $CO_2$ as fluid 2,
the Peng-Robinson bulk free energy density is defined by the following
\ben\bea{l}
h(\rho_1, \rho_2, T) =  \frac{r_m \rho_1 + \rho_2}{m_2} \varphi(T) -  \frac{r_m \rho_1 + \rho_2}{m_2} RT ln( \frac{m_2}{r_m \rho_1 + \rho_2} - b)
-  \\\\
\frac{r_m \rho_1 + \rho_2}{m_2}\frac{ a}{2\sqrt{2} b} ln[\frac{m_2 + (r_m \rho_1 + \rho_2) b(1+\sqrt{2})}{m_2 + (r_m \rho_1 + \rho_2) b(1-\sqrt{2})}]
+  \frac{r_m \rho_1 + \rho_2}{m_2} RT[\frac{r_m \rho_1}{r_m \rho_1 + \rho_2}ln\frac{r_m \rho_1}{r_m \rho_1 + \rho_2}+ \frac{\rho_2}{r_m \rho_1 + \rho_2}ln\frac{\rho_2}{r_m \rho_1 + \rho_2}].
\eea\label{PR_free_energy}
\een
where $R$ is the ideal gas constant, $\varphi(T) = -RT(1-log(\lambda^3))$ is a temperature-dependent function, $\lambda$ is the thermal wavelength of a massive particle, $m_i$ is the molar mass of component i for $i=1,2,$ respectively, $r_m = {m_2}/{m_1}$ is the ratio of the molar mass of carbon dioxide $m_2$ to the molar mass of n-decane $m_1$, $b(\rho_1, \rho_2)$ is a volume parameter and $a(\rho_1, \rho_2, T)$ is an interaction parameter. This free energy was proposed to extend that of the Van der Waals' to describe the deviation away from the ideal gas model. 

Another example of the bulk free energy density for polymeric liquids is given by the Flory-Huggins mixing energy density
\ben
h(\rho_1, \rho_2, T) =\frac{ k_B T}{m} [\frac{\rho_1}{N_1 } ln \frac{\rho_1}{\rho}+\frac{\rho_2}{N_2 } ln \frac{\rho_2}{\rho}+\chi \frac{\rho_1\rho_2}{\rho}],
\een
where $m$ is the mass of an average molecule in the mixture and $N_{1,2}$ are two polymerization indices.

Notice that $j_i, i=1,2$ in \eqref{compressible_sun} are obtained from the constitutive equation and if $\sum_{i=1}^2 j_i\neq 0$, this model  does not necessarily conserve mass locally. However, $\int_{V} (\rho_1+\rho_2) d{\bf x}$ is a constant. So, the  mass of the system is conserved globally.
This model describes a binary viscous compressible fluid systems in which mass is  conserved globally but not locally. In this model, the exact physical meaning of the velocity is lost due to the lack of local mass conservation. It is no long a mass average velocity! Therefore, what does the momentum equation stands for becomes fuzzy physically. The applicability of this model needs to be scrutinized further. A more general model can be built from \eqref{general} by specifying a more general mobility operator $\mathcal M$. However, we will not pursue it in this study.

Next, we impose the local mass conservation constraint to arrive at the model that conserves mass locally.
\subsection{Compressible model with local mass conservation law}
If  $j_1 + j_2 = 0$, the total mass of the system is conserved locally, i.e.,
\ben\bea{l}
\frac{\partial \rho}{\partial t} + \nabla \cdot (\rho {\bf v}) = 0,
\eea \label{Mass-conservation}
\een
which imposes an constraint on the mass fluxes:
\ben
\sum_{i=1}^2 \sum_{j=1}^2    \nabla \cdot  M_{ij} \cdot \nabla  \mu_j  = 0.
\een
We obtain the governing system of equations for the compressible fluid mixture as follows
\ben\bea{l}
\begin{cases}
\frac{\partial \rho}{\partial t} + \nabla \cdot (\rho {\bf v}) =0,\\\\
\frac{\partial \rho_i}{\partial t} + \nabla \cdot (\rho_i {\bf v}) =   \nabla \cdot  M_{i1} \cdot \nabla  \mu_1 +   \nabla \cdot  M_{i2} \cdot \nabla  \mu_2 , \quad i=1,2,\\\\
\frac{\partial (\rho {\bf v})}{\partial t}  +  \nabla \cdot (\rho {\bf v} {\bf v})  = 2 \nabla \cdot ( \eta{\bf D}) + \nabla (  {\nu} \nabla \cdot {\bf v})- \sum_{i=1}^{2} \rho_i \nabla  \mu_i.
\end{cases}
\eea \label{Pdes}
\een

Notice that we could have used $\rho$, $\rho_1$ as the fundamental variables in the derivation of the thermodynamic model in lieu of $\rho_1$ and $\rho_2$ since $\rho=\rho_1+\rho_2$.  With these variables, we reformulate the free energy density function
\ben\bea{l}
f(\rho_1, \rho_2, \nabla \rho_1, \nabla \rho_2) = f(\rho_1, \rho - \rho_1, \nabla \rho_1, \nabla (\rho - \rho_1))\\
\\
=\tilde f(\rho_1, \rho, \nabla \rho_1, \nabla \rho) = \tilde h(\rho_1, \rho, T)  + \frac{1}{2}(\tilde \kappa_{\rho_1 \rho_1} (\nabla \rho_1)^2 + 2\tilde \kappa_{\rho \rho_1} (\nabla \rho, \nabla \rho_1) + \tilde \kappa_{\rho \rho} (\nabla \rho)^2),
\eea\label{free_energy_local}
\een
where $\tilde \kappa_{\rho_1 \rho_1} = \kappa_{\rho_1 \rho_1} + \kappa_{\rho_2 \rho_2} - 2\kappa_{\rho_1 \rho_2}, \tilde \kappa_{\rho \rho_1} =  \kappa_{\rho_1 \rho_2} -  \kappa_{\rho_2 \rho_2}$, and $ \tilde \kappa_{\rho \rho} = \kappa_{\rho_2 \rho_2}$, where $\kappa_{\rho_1\rho_1}, \kappa_{\rho_1 \rho_2}, \kappa_{\rho_2 \rho_2}$ are the coefficients of the gradient terms in free energy \eqref{eq:kappa_coef}. 
The corresponding chemical potentials are given by
\ben
\bea{l}
{\tilde \mu_1}=\frac{\delta \tilde f}{\delta \rho_1}=\frac{\delta f}{\delta \rho_1} + \frac{\delta f}{\delta \rho_2} \frac{\partial \rho_2}{\partial \rho_1} = \mu_1 -  \mu_2,  \qquad
{\tilde \mu}=\frac{\delta \tilde f}{\delta \rho}=\frac{\delta f}{\delta \rho_2} \frac{\partial \rho_2}{\partial \rho} = \mu_2.\\
\mu_1={\tilde \mu}_1 +{\tilde \mu}, \qquad  \mu_2={\tilde \mu}.
\eea
\een
System (\eqref{Pdes}) reduces to
\ben\bea{l}
\begin{cases}
\frac{\partial \rho}{\partial t} + \nabla \cdot(\rho {\bf v})  = 0,\\\\
\frac{\partial \rho_1}{\partial t} + \nabla \cdot (\rho_1 {\bf v}) =   \nabla \cdot [M_{11} \cdot \nabla \tilde \mu_1+(M_{11}+M_{12})\cdot \nabla \tilde \mu],\\\\
\frac{\partial (\rho {\bf v})}{\partial t}  +  \nabla \cdot (\rho {\bf v} {\bf v})  = 2 \nabla \cdot ( \eta{\bf D}) + \nabla ( {\nu} \nabla \cdot {\bf v})-  \rho_1 \nabla \tilde \mu_1 - \rho \nabla \tilde  \mu.
\end{cases}
\eea \label{Pdes11}
\een
If we assign
\ben\bea{l}
 M_{12} =   M_{21} = -  M_{11}, \quad  M_{22} =  M_{11},
\eea\een
 system (\ref{Pdes}) reduces further to a special model
\ben\bea{l}
\begin{cases}
\frac{\partial \rho}{\partial t} + \nabla \cdot(\rho {\bf v})  = 0,\\\\
\frac{\partial \rho_1}{\partial t} + \nabla \cdot (\rho_1 {\bf v}) =   \nabla \cdot M_{11} \cdot \nabla \tilde \mu_1,\\\\
\frac{\partial (\rho {\bf v})}{\partial t}  +  \nabla \cdot (\rho {\bf v} {\bf v})  = 2 \nabla \cdot ( \eta{\bf D}) + \nabla ( {\nu} \nabla \cdot {\bf v})-  \rho_1 \nabla \tilde \mu_1 - \rho \nabla \tilde  \mu.
\end{cases}
\eea \label{Pdes2}
\een
This is a special model for compressible binary fluid mixtures among infinitely many choices in the mobility matrix. Apparently, model \eqref{Pdes11} is more general.

The boundary conditions at a solid boundary are given by (\ref{pdes-bc}) except that the last one is replaced by ${\bf n}\cdot \frac{\partial f}{\partial \nabla \rho}=0$ equivalently when $\rho$ is used as a fundamental variable. The energy dissipation rate of the special model  reduces to
\ben\bea{l}
\frac{dE_{total}}{dt}
=  -\int_{V}[ 2\eta {\bf D} :{\bf D} +   {\nu}  tr({\bf D})^2 +  \tilde \mu_1  M_{11} \tilde \mu_1
]  d{\bf x}\leq 0,
\eea\een
provided $\eta, {\nu} \geq 0$ and $M_{11} > 0$. This is a compressible binary fluid model that respects mass and momentum conservation.
For the more general model \eqref{Pdes11}, the energy dissipation property is warranted so long as the mobility matrix $\bf M$ is non-negative definite. So, this class of models is thermodynamically consistent.

We next show how this (special) model reduces to another class of compressible models when the two fluid components are incompressible, known as the quasi-incompressible model \cite{LowengrubRSA1998, Li&WangJAM2014}. For the more general compressible with a local mass conservation law, an analogous result can be obtained.

\subsection{Quasi-incompressible model}

When the fluid mixture is consisted of two incompressible viscous fluid components, where the specific densities $\hat{\rho_1}$ and $\hat{\rho_2}$ are constants, we  denote the volume fraction of fluid component 1 as $\phi$ and the other by $1-\phi$.
Then, the densities of the two fluids in the mixture are given as follows
\ben
\rho_1=\phi \hat{\rho}_1, \quad \rho_2=(1-\phi)\hat{\rho}_2.
\een
The total density of the fluid mixture is given by
\ben
\rho=\phi \hat{\rho}_1+(1-\phi)\hat{\rho}_2.
\een
If we use $\rho_1$ as a fundamental physical variable, $\rho$ is represented by $\rho_1$ as follows,
\ben
\rho=\rho_1+(1-\frac{\rho_1}{\hat{\rho}_1})\hat{\rho_2}=\hat{\rho}_2+(1-\frac{\hat{\rho}_2}{\hat{\rho}_1})\rho_1.\label{const}
\een
This means that the two variables $\rho$ and $\rho_1$ are related linearly in this fluid mixture system.
We view this as a special case of the fully compressible model subject to constraint (\ref{const}). To accommodate the constraint, we augment the free energy density by $\pi(\hat{\rho}_2+(1-\frac{\hat{\rho}_2}{\hat{\rho}_1})\rho_1-\rho)$, where $\pi$ is a Lagrange multiplier. We denote the modified free energy density function as $\hat{f}$,
\ben\bea{l}
\hat{f} = \tilde{f} (\rho_1, \nabla \rho_1, \rho, \nabla \rho)+ \pi [\hat{\rho}_2+(1-\frac{\hat{\rho}_2}{\hat{\rho}_1})\rho_1-\rho].
\eea\een

The corresponding chemical potentials and their relations to the chemical potentials in the  compressible model are given as follows
\ben\bea{l}
\hat{\mu}_1 = \frac{\delta \hat{f}}{\delta \rho_1} = \tilde \mu_1 + \pi(1-\frac{\hat{\rho}_2}{\hat{\rho}_1}),\quad
\hat{\mu} = \frac{\delta \hat{f}}{\delta \rho} =\tilde \mu -  \pi, \quad
\tilde \mu_1 = \frac{\delta  \tilde{f} }{\delta \rho_1}|_{\rho}, \\\\
\tilde \mu = \frac{\delta  \tilde{f} }{\delta \rho}|_{\rho_1},\quad
{\mu_{\phi}} =  \frac{\delta \hat{f}}{\delta \rho_1}|_{\rho} \frac{\delta \rho_1 }{\delta \phi} +  \frac{\delta \hat{f}}{\delta \rho}|_{\rho_1} \frac{\delta \rho }{\delta \phi}=  \hat{\rho}_1 \hat{\mu}_1 + (\hat{\rho}_1 - \hat{\rho}_2) \hat{\mu}.
\eea\een
From the mass conservation of the mixture system \eqref{Pdes2}-1, we have
\ben\bea{l}
(\hat{\rho}_1-\hat{\rho}_2) [\frac{\partial \phi}{\partial t}+\nabla \cdot (\phi {\bf v})]+\hat{\rho}_2 \nabla \cdot {\bf v}=0.
\eea\label{mass}
\een
The transport equation of $\rho_1$ is rewritten into
\ben\bea{l}
 \frac{\partial \phi}{\partial t} + \nabla \cdot (\phi {\bf v}) = \frac{1}{\hat{\rho}_1 } (\nabla \cdot M_{11} \cdot \nabla) (\hat {\mu_1} ).
\eea \label{Transport_phi}
\een
The linear momentum conservation equation is rewritten into
\ben
\rho ( \frac{\partial{\bf v}}{\partial t} +  {\bf v}\cdot \nabla {\bf v})= \nabla \cdot (2 \eta {\bf D}) + \nabla (  {\nu} \nabla \cdot {\bf v}) - \nabla \Pi-\phi \nabla \mu_{\phi},
\een
where $\eta$, ${\nu}$ are volume averaged viscosity coefficients and the  hydrostatic pressure is defined by
\ben\bea{c}
\Pi =  \hat{\rho_2}(\tilde \mu-\pi).
\eea\een
With this definition, the transport equation \eqref{Transport_phi} for $\phi$ is written into
\ben\bea{l}
 \frac{\partial \phi}{\partial t} + \nabla \cdot (\phi {\bf v}) = \frac{1}{\hat{\rho}_1^2 }(\nabla \cdot M_{11}\cdot  \nabla) (\mu_{\phi}+\Pi(1-\frac{\hat{\rho}_1}{\hat{\rho}_2})).
\eea \label{transport_2}
\een
Combining the mass conservation law \eqref{mass} and transport equation \eqref{transport_2} of the $\phi$, we obtain
\ben\bea{l}
\nabla \cdot {\bf v} = (1 - \frac{\hat{\rho}_1}{\hat{\rho}_2}) \frac{1}{\hat{\rho}_1^2} (\nabla \cdot M_{11} \cdot \nabla) (\mu_{\phi}+\Pi(1-\frac{\hat{\rho}_1}{\hat{\rho}_2})).
\eea\een

We summarize the governing equations of the quasi-incompressible model as follows
\ben\bea{l}
\begin{cases}
\nabla \cdot {\bf v} = (1 - \frac{\hat{\rho}_1}{\hat{\rho}_2}) \frac{1}{\hat{\rho}_1^2} (\nabla \cdot M_{11} \cdot \nabla )( \mu_{\phi}+\Pi(1-\frac{\hat{\rho}_1}{\hat{\rho}_2})),\\
\\
 \frac{\partial \phi}{\partial t} +  \nabla \cdot (\phi {\bf v}) = \frac{1}{\hat{\rho}_1^2}(\nabla \cdot M_{11} \cdot \nabla) (  \mu_{\phi} + \Pi(1-\frac{\hat{\rho}_1}{\hat{\rho}_2})),\\
\\
\rho  [\frac{\partial {\bf v}}{\partial t} + {\bf v}\cdot \nabla {\bf v}]= \nabla \cdot (2 \eta {\bf D}) + \nabla ( {\nu} \nabla \cdot {\bf v}) - \nabla \Pi-\phi \nabla \mu_{\phi} .
\end{cases}
\eea \label{Two_quasi1}
\een
The free energy density reduces to
\ben\bea{l}
\tilde{f} (\rho_1, \rho , \nabla \rho_1, \nabla \rho) = \tilde {h}(\hat{\rho}_1 \phi, \hat{\rho}_1 \phi+\hat{\rho}_2 (1-\phi), T) \\
+ \frac{1}{2}(\tilde \kappa_{\rho_1 \rho_1} (\nabla \rho_1)^2 + 2\tilde \kappa_{\rho_1 \rho} (\nabla \rho_1, \nabla \rho) + \tilde \kappa_{\rho \rho} (\nabla \rho)^2)\\\\
= \hat{h}(\phi)+\frac{1}{2} \hat{ \kappa}_{\phi \phi} \|\nabla \phi\|^2,
\eea \label{eq:kappa_coefqs}
\een
where $\hat{h}(\phi)=\tilde {h}(\hat{\rho}_1 \phi, (\hat{\rho}_1-\hat{\rho}_2) \phi+\hat{\rho}_2, T), \hat{\kappa}_{\phi \phi}=\tilde \kappa_{\rho_1 \rho_1} \hat{\rho_1}^2+2\tilde \kappa_{\rho_1 \rho} \hat{\rho_1} (\hat{\rho_1}-\hat{\rho_2})+\tilde \kappa_{\rho \rho} (\hat{\rho_1}-\hat{\rho_2})^2.$
This is the equation system for quasi-incompressible binary fluids obtained in \cite{Li&WangJAM2014}. The upshot of the derivation shows that we can obtain the constrained theory from the unconstrained theory by augmenting the free energy with the algebraic constraint via a Lagrange multiplier.

The energy dissipation rate of the binary quasi-incompressible fluid flow \eqref{Two_quasi1} is given by
\ben\bea{l}
\frac{dE_{total}}{dt}
=  -\int_{V}[ 2\eta {\bf D} :{\bf D} + {\nu} tr({\bf D})^2 +  \nabla \hat{\mu}_1 \cdot M_{11} \cdot \nabla \hat{\mu}_1   ]  d{\bf x}\leq 0,
\eea\een
provided $\eta, {\nu} \geq 0 $, $M_{11} > 0$, where $\hat{\mu}_1 = \frac{1}{\hat{\rho_1}} (\mu_{\phi} + \Pi(1-\frac{\hat{\rho}_1}{\hat{\rho}_2})   )$.

When $\hat{\rho}_1 = \hat{\rho}_2 = \rho$, the system reduces to an incompressible model
\ben\bea{l}
\begin{cases}
 \nabla \cdot {\bf v}=0,\\
\\
 \frac{\partial \phi}{\partial t} +  \nabla \cdot (\phi {\bf v}) = \frac{1}{{\rho}^2}(\nabla \cdot M_{11} \cdot \nabla)   \mu_{\phi},\\
\\
\rho  [\frac{\partial{\bf v}}{\partial t} + {\bf v}\cdot \nabla {\bf v}]= \nabla \cdot (2 \eta {\bf D}) + \nabla ( {\nu} \nabla \cdot {\bf v}) -  \nabla \Pi-\phi \nabla \mu_{\phi} .
\end{cases}
\eea\een
This is the incompressible model derived by Halperin et al \cite{Hohenberg&Halperin1977}.

These derivations can be readily extended to account for multi-component fluid systems.

\section{Hydrodynamic phase Field Models for N-component Multiphase Compressible Fluid Flows }

When   fluid mixtures are composed of   N fluid components, we use $\rho_i, i=1,2,\cdots, N$ to denote the mass density of the ith  component and assume the free energy of the fluid mixture is given by
\ben
F=\int_{\Omega} f(\rho_1, \nabla \rho_1, \cdots, \rho_N, \nabla \rho_N) d\bf x,
\een
where $f$ is the free energy density. The derivation of the hydrodynamic phase field models follows the procedures alluded to in the previous section. We  present the results next.


\subsection{Compressible model with the global mass conservation law}

We choose $\rho_1, \cdots, \rho_{N}$ as the primitive variables. Following the procedure outlined in the previous section, we obtain the governing system of equations for the N-component multi-phase viscous fluid mixture as follows
\ben\bea{l}
\begin{cases}
\frac{\partial \rho_i}{\partial t} + \nabla \cdot (\rho_i {\bf v}) = j_i =  \sum_{j=1}^{N} \nabla \cdot M_{ij} \cdot \nabla \mu_j, \quad i=1,2,\cdots, N,\\\\
\frac{\partial (\rho {\bf v})}{\partial t}  +  \nabla \cdot (\rho {\bf v} {\bf v}) - \frac{1}{2}(\sum_{i=1}^{N}j_i){\bf v}  = 2 \nabla \cdot ( \eta{\bf D}) + \nabla (  {\nu} \nabla \cdot {\bf v})-  \sum_{i=1}^{N} \rho_i \nabla \mu_i,
\end{cases}
\eea
\een
where $M_{ij}$, i, j = 1, ..., N, are the mobility coefficients, and $\eta=\sum_{i=1}^{N} \eta_i \frac{\rho_i}{\rho}, {\nu}=\sum_{i=1}^{N} {\nu_i} \frac{\rho_i}{\rho}$ are mass-average viscosities, respectively.

The energy dissipation rate  is given by
\ben\bea{l}
\frac{dE_{total}}{dt}
=  -\int_{V}[ 2\eta {\bf D} :{\bf D} +   {\nu} tr({\bf D})^2 \\

+ (\nabla \mu_1,\nabla  \mu_2, \cdots,  \nabla \mu_{N}) \cdot   {\bf M} \cdot (\nabla  \mu_1, \nabla  \mu_2, \cdots, \nabla  \mu_{N})]  d{\bf x}\leq 0,
\eea\een
provided $\eta, {\nu} \geq 0$, $ {\bf M}=(M_{ij})_{ij=1}^N$ is a symmetric non-negative definite mobility coefficient matrix.

\subsection{Compressible model with the local mass conservation law}
If  $\sum_{i=1}^{N}j_i = 0$, the total mass of the system is conserved locally, i.e.,
\ben\bea{l}
\frac{\partial \rho}{ \partial t} + \nabla \cdot (\rho {\bf v}) = 0,
\eea \label{Mass-conservation_Multi}
\een
We obtain the governing system of equations as follows
\ben\bea{l}
\begin{cases}
\frac{\partial \rho}{\partial t} + \nabla \cdot (\rho {\bf v}) = 0, \quad \hbox{or} \quad \sum_{i=1}^N \sum_{j=1}^N   \nabla \cdot  M_{ij} \cdot \nabla  \mu_j  = 0,\\\\
\frac{\partial \rho_i}{\partial t} + \nabla \cdot (\rho_i {\bf v}) =   \sum_{j=1}^{N} \nabla \cdot M_{ij} \cdot \nabla \mu_j, \quad i=1,2,\cdots,N,\\\\
\frac{\partial (\rho {\bf v})}{\partial t}  +  \nabla \cdot (\rho {\bf v} {\bf v})  = 2 \nabla \cdot ( \eta{\bf D}) + \nabla (  {\nu} \nabla \cdot {\bf v})- \sum_{i=1}^{N} \rho_i \nabla  \mu_i,
\end{cases}
\eea \label{Pdes_Multi}
\een
where $\eta, {\nu}$ are mass averaged shear and volumetric viscosities, $\bf v$ is the mass average velocity and ${\bf M}=(M_{ij})_{i,j=1}^N$ is the symmetric mobility coefficient matrix.
In this case, the energy dissipation rate  is given by
\ben\bea{l}
\frac{dE_{total}}{dt}
=  -\int_{V}[ 2\eta {\bf D} :{\bf D} +{\nu} tr({\bf D})^2 \\

+ (\nabla  \mu_1,\nabla   \mu_2, \cdots,  \nabla \mu_{N}) \cdot  {\bf M} \cdot (\nabla  \mu_1, \nabla   \mu_2, \cdots, \nabla \mu_{N})]  d{\bf x}\leq 0,
\eea\een
provided $\eta, {\nu} \geq 0$ and ${\bf M}$ is a symmetric non-negative definite mobility coefficient matrix  subject to constraint $\quad \sum_{i=1}^N \sum_{j=1}^N   \nabla \cdot  M_{ij} \cdot \nabla  \mu_j  = 0$.

Analogously, we choose $\rho_1, \cdots, \rho_{N-1}, \rho$ as the primitive variables, where $\rho=\sum_{i=1}^N \rho_i.$ Then, we represent $\rho_N=\rho-\sum_{i=1}^{N-1} \rho_i$. The free energy density is written as
\ben\bea{l}
f(\rho_1, \nabla \rho_1, \cdots, \rho_N, \nabla \rho_N) = f(\rho_1, \nabla \rho_1, \cdots, \rho - \sum_{i=1}^{N-1} \rho_i, \nabla(\rho - \sum_{i=1}^{N-1} \rho_i))\\
\\
=\tilde f(\rho_1, \nabla \rho_1, \cdots, \rho_{N-1}, \nabla \rho),
\eea\een
The corresponding chemical potentials are given by
\ben
\bea{l}
\tilde \mu_i=\frac{\delta \tilde f}{\delta \rho_i}=\frac{\delta f}{\delta \rho_i}   +   \frac{\delta f}{\delta \rho_N}  \frac{\delta \rho_N}{\delta \rho_i} = \mu_i - \mu_N, i=1, \cdots, N-1, \tilde \mu=\frac{\delta \tilde f}{\delta \rho}=\frac{\delta f}{\delta \rho_N} \frac{\delta \rho_N}{\delta \rho}= \mu_N.\\
\mu_i=\tilde \mu_i+\tilde \mu, \mu_N=\tilde \mu.
\eea
\een
The transport equation of the densities are given by
\ben
\frac{\partial \rho_i}{\partial t} + \nabla \cdot (\rho_i {\bf v}) =   \sum_{j=1}^{N-1} \nabla \cdot  M_{ij} \cdot \nabla \tilde \mu_j+(\sum_{j=1}^N \nabla \cdot M_{ij}) \tilde \mu,  \quad i=1,2,\cdots,N-1.
\een
The  mass conservation equation implies
\ben
\sum_{i=1}^N\sum_{j=1}^{N-1} \nabla \cdot M_{ij}\cdot \nabla \tilde \mu_j+(\sum_{i,j=1}^N \nabla \cdot M_{ij}\cdot \nabla)\tilde \mu=0.
\een
The mobility coefficients must satisfy the above constraint.
If we assign
\ben\bea{l}
M_{iN} = - \sum_{j=1}^{N-1}   M_{ij} = M_{Ni}, \quad
M_{NN} = - \sum_{i=1}^{N-1}  M_{iN} = \sum_{i=1}^{N-1} \sum_{j=1}^{N-1} M_{ij},
\eea\een
the constraint is satisfied and  system (\ref{Pdes_Multi}) reduces to a special model
\ben\bea{l}
\begin{cases}
\frac{\partial \rho}{\partial t} + \nabla \cdot(\rho {\bf v})  = 0,\\\\
\frac{\partial \rho_i}{\partial t} + \nabla \cdot (\rho_i {\bf v}) =   \sum_{j=1}^{N-1} \nabla \cdot  M_{ij} \cdot \nabla \tilde \mu_j,  \quad i=1,2,\cdots,N-1,\\\\
\frac{\partial (\rho {\bf v})}{\partial t}  +  \nabla \cdot (\rho {\bf v} {\bf v})  = 2 \nabla \cdot ( \eta{\bf D}) + \nabla ( {\nu} \nabla \cdot {\bf v})- \sum_{i=1}^{N-1} \rho_i \nabla \tilde \mu_i - \rho \nabla \tilde \mu.
\end{cases}
\eea \label{Multi-compressible}
\een
This is a special model for compressible fluid mixtures of N-components.
The energy dissipation rate  is given by
\ben\bea{l}
\frac{dE_{total}}{dt}
=  -\int_{V}[ 2\eta {\bf D} :{\bf D} +  {\nu} tr({\bf D})^2 \\

+ (\nabla \tilde \mu_1,\nabla \tilde \mu_2, \cdots, \nabla  \tilde \mu_{N-1}) \cdot {\bf M} \cdot (\nabla \tilde \mu_1, \nabla \tilde \mu_2, \cdots, \nabla \tilde \mu_{N-1})]  d{\bf x}\leq 0,
\eea\een
provided $\eta, {\nu} \geq 0$ and ${\bf M} = (M_{ij})_{i,j=1}^{N-1}$ is a symmetric non-negative definite mobility coefficient matrix. 


\subsection{Quasi-incompressible  model}

When each of the fluid component is incompressible in the viscous fluid mixture, we denote the volume fraction of the ith component as $\phi_{i}$ and specific density as $\hat{\rho_i}$ for $i=1, \cdots,N$, respectively. Then, $\sum_{i=1}^N\phi_i=1$ and the total mass density in the mixture is given by
\ben
\rho=\sum_{i=1}^N \phi_i \hat{\rho_i}=\sum_{i=1}^{N-1}{\phi_i} \hat{\rho_i}+(1-\sum_{i=1}^{N-1}{\phi_i}) \hat{\rho_N} = \sum_{i=1}^{N-1}{\rho_i}+(1-\sum_{i=1}^{N-1}{\frac{\rho_i}{\hat{\rho_i}}}) \hat{\rho_N}.
\een
We assume the volume fraction of the Noth component is nonzero. Then, the free energy density is a functional of the first $N-1$ volume fractions $(\phi_1, \cdots, \phi_{N-1})$. If we augment the free energy by $\pi(\sum_{i=1}^{N-1}{\rho_i}+(1-\sum_{i=1}^{N-1}{\frac{\rho_i}{\hat{\rho_i}}}) \hat{\rho_N}-\rho)$, where $\pi$ is a Lagrange multiplier, then, the modified free energy density function is given by
\ben\bea{l}
\hat{f} = \tilde{f}(\rho_1, \nabla \rho_1,  ..., \rho_{N-1}, \nabla \rho_{N-1}, \rho, \nabla \rho) + \pi [\sum_{i=1}^{N-1}{\rho_i}+(1-\sum_{i=1}^{N-1}{\frac{\rho_i}{\hat{\rho_i}}}) \hat{\rho_N}-\rho] .
\eea\een
Following the procedure alluded to in the previous section, we derive the following
governing system of equations of the quasi-incompressible fluid from the special compressible model as follows
\ben\bea{l}
\begin{cases}
\nabla \cdot {\bf v} = \sum_{i=1}^{N-1}  \sum_{j=1}^{N-1}(1 - \frac{\hat{\rho}_j}{\hat{\rho}_N}) \frac{1}{\hat{\rho}_i\hat{\rho}_j} (\nabla \cdot M_{ij} \cdot \nabla )( \mu_{\phi_j}+ \Pi (1-\frac{\hat{\rho}_j}{\hat{\rho}_N})),\\
\\
 \frac{\partial \phi_i}{\partial t} +  \nabla \cdot (\phi_i {\bf v}) =\sum_{j=1}^{N-1} \frac{1}{\hat{\rho}_i \hat{\rho}_j}(\nabla \cdot M_{ij} \cdot \nabla) ( \mu_{\phi_j} + \Pi(1-\frac{\hat{\rho}_j}{\hat{\rho}_N})),  i = 1, 2, \cdots, N-1,\\
\\
\rho  [\frac{\partial {\bf v}}{\partial t} + {\bf v}\cdot \nabla {\bf v}]= \nabla \cdot (2 \eta {\bf D}) + \nabla ( {\nu} \nabla \cdot {\bf v}) - \nabla \Pi - \sum_{i=1}^{N-1}\phi_i \nabla \mu_{\phi_i},
\end{cases}
\eea \label{Multi-quasi1}
\een
where
\ben\bea{l}
\tilde{\mu}_i =  \frac{\delta \tilde{f} }{\delta \rho_i},  \quad \tilde{\mu} = \frac{\delta  \tilde{f} }{\delta \rho},   \\

\hat{\mu}_i =  \tilde{\mu}_i + \pi(1-\frac{\hat{\rho}_N}{\hat{\rho}_i}), \quad i=1, \cdots, N-1,  \quad \hat{\mu} = \tilde{\mu} - \pi,  \\

\mu_{\phi_i} = \frac{\delta \hat{f}}{\delta \phi_i} =  \frac{\delta \hat{f}}{\delta \rho_i}|_{\rho} \frac{\delta \rho_i}{\delta \phi_i} + \frac{\delta \hat{f}}{\delta \rho}|_{\rho_i} \frac{\delta \rho}{\delta \phi_i}=  \hat{\rho}_i \hat{\mu}_i +  (\hat{\rho}_i - \hat{\rho}_N) \hat{\mu},\\

\Pi = - \hat{\rho}_N \pi + \hat{\rho}_N \mu,
\eea\een
and $\Pi$ serves as the hydrostatic pressure.

The energy dissipation rate is
\ben\bea{l}
\frac{dE_{total}}{dt}
=  -\int_{V}[ 2\eta {\bf D} :{\bf D} +{\nu} tr({\bf D})^2 \\

+ (\nabla \hat{\mu}_1, \nabla \hat{\mu}_2, \cdots, \nabla \hat{\mu}_{N-1}) \cdot {\bf M} \cdot (\nabla \hat{\mu}_1, \nabla \hat{\mu}_2, \cdots, \nabla \hat{\mu}_{N-1})

  ]  d{\bf x}\leq 0,
\eea\een
provided $\eta, {\nu} \geq 0$, ${\bf M} = (M)_{i,j=1}^{N-1}$ is a symmetric non-negative definite matrix, where
$\hat{\mu}_i = \frac{1}{\hat{\rho_i}}[ \mu_{\phi_j}+  \Pi (1-\frac{\hat{\rho}_j}{\hat{\rho}_N}) ]$.
A more general model can be derived from the general compressible model by enforcing the incompressibility  constraint. But, we will not present it here.

For a fluid mixture with $\hat{\rho}_i = \rho$, the system reduces to an incompressible model
\ben\bea{l}
\begin{cases}
 \nabla \cdot {\bf v}=0,\\
\\
 \frac{\partial \phi_i}{\partial t} +  \nabla \cdot (\phi_i {\bf v}) =\sum_{j=1}^{N-1} \frac{1}{{\rho}^2}(\nabla \cdot M_{ij} \cdot \nabla)   \mu_{\phi_j}, \quad i = 1, 2, \cdots, N-1,\\
\\
\rho  [\frac{\partial{\bf v}}{\partial t} + {\bf v}\cdot \nabla {\bf v}]= \nabla \cdot (2 \eta {\bf D}) + \nabla (  {\nu} \nabla \cdot {\bf v}) -  \nabla \Pi-   \sum_{i=1}^{N-1}\phi_i \nabla \mu_{\phi_i} .
\end{cases}
\eea\label{Multi-quasi2}
\een

For phase field models of N components where $N\geq 2$, there exists a second way to derive the quasi-incompressible phase field model. We begin with a fully compressible model of $N+1$ components, each of which is of density $\rho_i, i=1,\cdots,N$ and $\rho$.  We assume the free energy density depends on  $(\rho_1, \cdots, \rho_N, \rho)$.
 The second approach to derive the quasi-incompressible model is to augment the free energy by $\pi(\sum_{i=1}^{N}{\rho_i}-\rho) + B(\sum_{i=1}^{N}{\frac{\rho_i}{\hat{\rho_i}}}- 1)$, where $\pi$ and B are two Lagrange multipliers. We define the modified free energy density function by
\ben\bea{l}
\hat{f} = f(\rho_1, \nabla \rho_1,  \cdots, \rho_N, \nabla \rho_N) + \pi(\sum_{i=1}^{N}{\rho_i}-\rho) + B(\sum_{i=1}^{N}{\frac{\rho_i}{\hat{\rho_i}}}- 1).
\eea\label{free_energy_quasi2}
\een
The chemical potentials are given by
\ben\bea{l}
\hat{\mu}_i = \frac{\delta \hat{f}}{\delta \rho_i}  = \frac{\delta {f}}{\delta \rho_i} + \frac{1}{\hat{\rho_i}} B + \pi=  \mu_i + \frac{1}{\hat{\rho_i}} B + \pi, \quad i=1,\cdots, N,\\

\hat{\mu} =  \frac{\delta \hat{f}}{\delta \rho}  = - \pi.
\eea\een
The governing   system of equations with N+1 components subject to the two constraints is given by
\ben\bea{l}
\begin{cases}
\frac{\partial \rho}{\partial t} + \nabla \cdot (\rho {\bf v}) = 0,\\\\
\frac{\partial \rho_i}{\partial t} + \nabla \cdot (\rho_i {\bf v}) =   \sum_{j=1}^N \nabla \cdot  M_{ij} \cdot \nabla \hat{ \mu}_j, \quad i=1, 2,\cdots, N,\\\\
\frac{\partial (\rho {\bf v})}{\partial t}  +  \nabla \cdot (\rho {\bf v} {\bf v})  = 2 \nabla \cdot ( \eta{\bf D}) + \nabla ({\nu} \nabla \cdot {\bf v})- \sum_{i=1}^{N} \rho_i \nabla  \hat{\mu}_i -  \rho \nabla \hat{\mu},
\end{cases}
\eea \label{N1pdes}
\een
where ${\bf M}$ is the symmetric   mobility matrix, which satisfies $\sum_{i,j=1}^N \nabla \cdot M_{ij} \cdot \nabla \hat{\mu}_j=0$. This is a more general quasi-incompressible model.

In fact, if we assign $M_{Ni} = M_{iN} = - \sum_{j=1}^{N-1} M_{i j}$ and apply the constraints $\sum_{i=1}^N \rho_i = \rho$, $\sum_{i=1}^N \phi_i = 1$ and $\rho_i = \phi_i \hat{\rho_i}$, we obtain the chemical potential with respect to $\phi_i$, i=1,2, $\cdots$, N-1, in the quasi-incompressible limit,
\ben\bea{l}
\mu_{\phi_i} =  \frac{\delta \hat{f}}{\delta \phi_i} =  \frac{\delta \hat{f}}{\delta \rho_i} \frac{\partial \rho_i}{\partial \phi_i} + \frac{\delta \hat{f}}{\delta \rho} \frac{\partial \rho}{\partial \phi_i} +\frac{\delta \hat{f}}{\delta \rho_N} \frac{\partial \rho_N}{\partial \phi_i} = \hat{\mu_i} \hat{\rho}_i + \hat{\mu} (\hat{\rho_i} - \hat{\rho_N}) - \hat{\rho_N} \hat{\mu_N},  i=1,\cdots, N.
\eea\een
If we define
\ben\bea{l}
\Pi = \hat{\rho_N} (\hat{\mu} + \hat{\mu}_N) = \hat{\rho_N} {\mu}_N + B,
\eea\een
The model in (\ref{N1pdes}) reduces to
\ben\bea{l}
\begin{cases}
\nabla \cdot {\bf v} = \sum_{i=1}^{N-1}  \sum_{j=1}^{N-1}  (1 - \frac{\hat{\rho}_j}{\hat{\rho}_N}) \frac{1}{\hat{\rho}_i \hat{\rho}_j} (\nabla \cdot M_{ij} \cdot \nabla )( \mu_{\phi_j}+  \Pi (1-\frac{\hat{\rho}_j}{\hat{\rho}_N})),\\\\

 \frac{\partial \phi_i}{\partial t} +  \nabla \cdot (\phi_i {\bf v}) =\sum_{j=1}^{N-1} \frac{1}{\hat{\rho}_i \hat{\rho}_j}(\nabla \cdot M_{ij} \cdot \nabla) (  \mu_{\phi_j} +   \Pi(1-\frac{\hat{\rho}_j}{\hat{\rho}_N})),  i = 1, 2, \cdots, N-1,\\
\\
\rho  [\frac{\partial {\bf v}}{\partial t} + {\bf v}\cdot \nabla {\bf v}]= \nabla \cdot (2 \eta {\bf D}) + \nabla ( {\nu} \nabla \cdot {\bf v}) - \nabla   \Pi - \sum_{i=1}^{N-1}\phi_i \nabla \mu_{\phi_i},
\end{cases}
\eea
\een
which is exactly the quasi-incompressible model given in (\ref{Multi-quasi1}).


\section{Non-dimensionalization}

Next, we non-dimensionalize the binary model equations and compare their near equilibrium dynamics.

\subsection{Compressible model with the global mass conservation law}

In model \eqref{compressible_sun}, selecting characteristic time scale $t_0$, characteristic length scale $l_0$, and characteristic density scale $\rho_0$, we nondimensionalize the variables and parameters as follows
\ben\bea{l}
\tilde{t} = \frac{t}{t_0}, \quad \tilde{x} = \frac{x}{l_0}, \quad \tilde{\rho}_1 = \frac{\rho_1}{\rho_0}, \quad \tilde{\rho}_2 = \frac{\rho_2}{\rho_0} , \quad \tilde{{\bf v}} = \frac{{\bf v}t_0}{l_0}, \quad \tilde{M}_{ij}  = \frac{M_{ij} }{t_0 \rho_0}, \quad i, j = 1, 2, \\

\frac{1}{{Re}_s} = \tilde{\eta} = \frac{t_0}{\rho_0 l_0^2} \eta,  \quad \frac{1}{{Re}_v} = \tilde{{\nu}} = \frac{t_0}{\rho_0 l_0^2} {\nu}, \quad \tilde{\mu}_1 = \frac{t_0^2}{ l_0^2} \mu_1, \quad \tilde{\mu}_2 = \frac{t_0^2}{ l_0^2} \mu_2, \quad J_i=\frac{j_i t_0}{\rho_0}, i=1,2,
\eea\een
where $Re_s$, $Re_v$ are the Reynolds number corresponding to the shear and volumetric stresses. The scaling of chemical potentials $\mu_1$, $\mu$ results from the non-dimensionalization of the total energy.
We summarize the governing equation with non-dimensional variables and parameters as follows, dropping the $\tilde{}$ for simplicity,
\ben\bea{l}
\begin{cases}
\frac{\partial \rho_1}{\partial t} + \nabla \cdot (\rho_1 {\bf v}) = J_1 =  \nabla \cdot M_{11} \cdot \nabla \mu_1 +  \nabla \cdot M_{12} \cdot \nabla \mu_2 ,\\
\\
\frac{\partial \rho_2}{\partial t} + \nabla \cdot (\rho_2 {\bf v}) = J_2 =  \nabla \cdot M_{12} \cdot \nabla \mu_1 + \nabla \cdot M_{22} \cdot \nabla \mu_2,\\
\\
\frac{\partial (\rho {\bf v})}{\partial t} + \nabla \cdot (\rho {\bf vv}) - \frac{1}{2}(J_1 + J_2) {\bf v}=
2 \nabla \cdot (  \frac{1}{{Re}_s}  {\bf D}) + \nabla (\frac{1}{{Re}_v} \nabla \cdot {\bf v}) - \rho_1 \nabla \mu_1 - \rho_2 \nabla \mu_2.
\end{cases}
\eea \label{eq:Nondim_com2}
\een
\subsection{Compressible model with the local mass conservation law}
Analogously, in model \eqref{Pdes2}, we nondimensionalize the variables and parameters as above and in particular
\ben\bea{l}
\quad \tilde{M}_{11}  = \frac{M_{11} }{t_0 \rho_0}.
\eea\een
We summarize the governing equation with non-dimensional variables and parameters as follows, dropping the $\tilde{}$ for simplicity,
\ben\bea{l}
\begin{cases}
\frac{\partial \rho}{\partial t} + \nabla \cdot (\rho {\bf v}) = 0,\\
\\
\frac{\partial \rho_1}{\partial t} + \nabla \cdot (\rho_1 {\bf v}) =  \nabla \cdot M_{11} \cdot \nabla \tilde{\mu}_1,\\
\\
\frac{\partial (\rho {\bf v})}{\partial t} + \nabla \cdot (\rho {\bf vv}) = 2 \nabla \cdot (\frac{1}{{Re}_s} {\bf D}) + \nabla (  \frac{1}{{Re}_v} \nabla \cdot {\bf v}) - \rho_1 \nabla \tilde{ \mu}_1 - \rho \nabla \tilde{\mu}.
\end{cases}
\eea \label{eq:Nondim_com1}
\een

\subsection{Quasi-incompressible model}
In model \eqref{Two_quasi1}, in addition to the above, we nondimensionalize two  new ones as follows:
\ben\bea{l}
\quad \tilde{\mu}_{\phi} = \frac{t_0^2}{\rho_0 l_0^2} \mu_{\phi}, \quad \tilde{\Pi} = \Pi \frac{t_0^2}{\rho_0 l_0^2}.
\eea\een
Dropping the $\tilde{}$ on the non-dimensionalized variables and parameters, the governing equation system of the quasi-incompressible fluid flows is written as follows,
\ben\bea{l}
\begin{cases}
\nabla \cdot {\bf v} = (1 - \frac{\hat{\rho}_1}{\hat{\rho}_2}) \frac{1}{\hat{\rho}_1^2} (\nabla \cdot M_{11} \cdot \nabla )( \mu_{\phi}+\Pi(1-\frac{\hat{\rho}_1}{\hat{\rho}_2})),\\
\\
 \frac{\partial \phi}{\partial t} +  \nabla \cdot (\phi {\bf v}) = \frac{1}{\hat{\rho}_1^2}(\nabla \cdot M_{11} \cdot \nabla) (  \mu_{\phi} + \Pi(1-\frac{\hat{\rho}_1}{\hat{\rho}_2})),\\
\\
\rho  [\frac{\partial {\bf v}}{\partial t} + {\bf v}\cdot \nabla {\bf v}]= 2  \nabla \cdot ( \frac{1}{{Re}_s} {\bf D}) + \nabla (  \frac{1}{{Re}_v} \nabla \cdot {\bf v}) - \nabla \Pi-\phi \nabla \mu_{\phi} .
\end{cases}
\eea\een


\section{Comparison of the models}

We investigate near equilibrium dynamics by conducting a linear stability analysis of the models from each class about a constant steady state. Through analyzing the dispersion relations of the selected models, we would like to identify  the intrinsic relation among compressible, quasi-incompressible and incompressible models, in particular, to reveal the consequence of the hierarchical reduction to linear stability. We focus on models of a binary fluid mixture only in this study.

\subsection{Linear stability analysis of the compressible model with the global mass conservation law}

This compressible model admits one constant solution:
\ben
{\bf v}={\bf 0},\quad \rho_1=\rho_1^0,\quad \rho_2=\rho_2^0,\label{csol}
\een
where $\rho_1^0, \rho_2^0$ are constants. We perturb the constant solution with the normal mode as follows:
\ben
{\bf v}= \epsilon e^{\alpha t + i {\bf k} \cdot {\bf x}} {\bf v}^{c}, \quad \rho_1=\rho_1^0 + \epsilon e^{\alpha t + i {\bf k} \cdot {\bf x}} {\rho_1}^{c},\quad  \rho_2=\rho_2^0 +\epsilon e^{\alpha t + i {\bf k} \cdot {\bf x}} {\rho_2}^{c} .
\een
where $\epsilon$ is a small parameter, representing the magnitude of the perturbation, and ${\bf v}^c, \rho_1^c, \rho_2^c$ are constants, $\alpha $ is the growth rate, $\bf k$ is the wave number of the  perturbation. Without loss of generality, we limit our study to 1 dimensional perturbations in $\bf k$ in 2D models. Substituting these perturbations into the equations in (\ref{eq:Nondim_com2}) and truncate the equations  at order $O(\epsilon)$, we obtain the linearized equations.  The dispersion equation of the linearized equation systems is given by the algebraic equation of $\alpha$:
\ben\bea{l}
(\frac{1}{{Re}_s} k^2 + \alpha \rho^0)\{\alpha^3 \rho_0 + \alpha^2 k^2[\frac{1}{{Re}} + \rho^0 M_{11} (h_{\rho_1 \rho_1} + \kappa_{\rho_1 \rho_1}k^2) +  \rho^0 M_{22} (h_{\rho_2 \rho_2} + \kappa_{\rho_2 \rho_2}k^2) ] \\

+ \alpha^2 k^2[ 2 \rho^0 M_{12}(h_{\rho_1 \rho_2} + \kappa_{\rho_1 \rho_2}k^2)] + \alpha [{\bf p}^T \cdot {\bf C} \cdot {\bf p} + {\bf p}^T\cdot  {\bf K} \cdot {\bf p}k^2 ]k^2 \\

+\alpha \frac{1}{{Re}} [  M_{11}  (h_{\rho_1 \rho_1} + \kappa_{\rho_1 \rho_1}k^2)  + M_{22} (h_{\rho_2 \rho_2} + \kappa_{\rho_2 \rho_2}k^2) + 2 M_{12}(h_{\rho_1 \rho_2} + \kappa_{ \rho_1 \rho_2}k^2) ]k^4
\\
+ \alpha \rho^0  |{\bf M}| [ (h_{\rho_1\rho_1} + \kappa_{\rho_1 \rho_1}k^2) (h_{\rho_2\rho_2} + \kappa_{\rho_2 \rho_2}k^2) -  (h_{\rho_1\rho_2} + \kappa_{\rho_1 \rho_2}k^2)^2]k^4 \\

+ k^4 ( \frac{1}{{Re}} |{\bf M}| k^2 + M_{22} (\rho_1^0)^2 + M_{11} (\rho_2^0)^2 - 2 M_{12} \rho_1^0 \rho_2^0  )\\

 [ (h_{\rho_1\rho_1} + \kappa_{\rho_1 \rho_1}k^2) (h_{\rho_2\rho_2} + \kappa_{\rho_2 \rho_2}k^2) -  (h_{\rho_1\rho_2} + \kappa_{\rho_1 \rho_2}k^2)^2]  \}
= 0.
\eea \label{eq:Algebraic_com1}
\een
where ${\bf p} = (\rho_1^0, \rho_2^0)^T$ and $\frac{1}{{Re}} = 2 \frac{1}{{Re}_s} + \frac{1}{{Re}_v}$. $|{\bf M}|$ is the determinant of the mobility coefficient matrix ${\bf M} = (M_{i,j})$
, $\bf K$ is the coefficient matrix of the conformational entropy
\ben
{\bf K} = \left(
\bea{cc}
\kappa_{\rho_1 \rho_1 }  & \kappa_{\rho_1 \rho_2}  \\
\kappa_{\rho_1 \rho_2}  & \kappa_{\rho_2 \rho_2 }
\eea
\right),
\een
$\bf C$ is the Hessian of the bulk free energy $h(\rho_1, \rho_2, T)$ in \eqref{eq:kappa_coef} with respect to  $\rho_1$ and $\rho_2$,
 \ben
{\bf C} = \left(
\bea{cc}
h_{\rho_1 \rho_1}  & h_{\rho_1 \rho_2}  \\
h_{\rho_1 \rho_2}  & h_{\rho_2 \rho_2}  \\
\eea
\right),
\een
where $h_{\rho_i \rho_j}$ represents the second order derivative of the bulk free energy density $h(\rho_1, \rho_2, T)$ with respect to $\rho_i$ and $\rho_j$, i, j = 1, 2.

One root of equation \eqref{eq:Algebraic_com1} is given by
\ben
\alpha_0=-\frac{1}{\rho^0 Re_s}k^2.
\een
This is the viscous mode associated to the viscous stress.
The other three roots are governed by a cubic polynomial equation and their closed forms are essentially impenetrable.
Instead,  we present them using asymptotic formulae in long and short wave range and numerical calculations in the intermediate wave range.

The asymptotic expressions of the three growth rates at  $|k| \ll 1$ are given by
\ben\bea{l}
\alpha_{1} = x_1 k^2 + y_1 k^4 + O(k^5),\quad
\alpha_{2,3} = x_{2,3} k + y_{2,3}k^2 + O(k^3),
\eea\een
where
\ben\bea{l}
x_1 = - \frac{g_1 |{\bf C}|}{{\bf p}^T \cdot {\bf C} \cdot {\bf p}},\quad x_{2,3} = \pm \sqrt{- \frac{{\bf p}^T \cdot {\bf C} \cdot {\bf p} }{\rho^0}}, \\
\\
y_1 = -\frac{1}{{\bf p}^T \cdot {\bf C} \cdot {\bf p}}\big[  \frac{1}{{Re}} |{\bf M}| |{\bf C}| + d g_1    \big]
-\frac{1}{{\bf p}^T \cdot {\bf C} \cdot {\bf p}}\big[ \rho^0 x_1^3 + x_1^2(\frac{1}{{Re}} + \rho^0 {\bf M} : {\bf C})  +\\
 x_1(\rho^0 |{\bf M}| |{\bf C}| + \frac{1}{{Re}} {\bf M} : {\bf C} + {\bf p}^T \cdot {\bf K} \cdot {\bf p}) \big], \\
\\
y_{2,3} = - \frac{1}{2 \rho^0 Re} - \frac{1}{2{\bf p}^T\cdot  {\bf C} \cdot {\bf p}}[ M_{11}(\rho_1^0 h_{\rho_1 \rho_1} + \rho_2^0 h_{\rho_1 \rho_2})^2 + M_{22} (\rho_1^0 h_{\rho_1 \rho_2} + \rho_2^0 h_{\rho_2 \rho_2})^2].
\eea\een
where $d = h_{\rho_1 \rho_1} \kappa_{\rho_2 \rho_2} + h_{\rho_2 \rho_2} \kappa_{\rho_1 \rho_1} - 2h_{\rho_1 \rho_2} \kappa_{\rho_1 \rho_2}$ and $g_1 = M_{22} (\rho_1^0 )^2 + M_{11} ( \rho_2^0)^2 - 2 M_{12}\rho_1^0 \rho_2^0 \geq 0$, since ${\bf M} \geq 0$ and at least one of its eigenvalues is positive.

When $|k| \gg 1$, the  three growth rates  are given by
\ben\bea{l}
\alpha_{1,2} = x_{1,2}k^4 + y_{1,2}k^2 + O(k),\quad
\alpha_3 = x_3 k^2 + y_3 + O(\frac{1}{k}),
\eea \label{eq:alpha_k_Inf_com1}
\een
where
\ben\bea{l}
x_{1,2} = - \frac{ {\bf M} : {\bf K} }{2} \pm \frac{1}{2\rho^0} \sqrt{({\bf M} :{\bf K} \rho^0)^2 - 4[ \frac{1}{{Re}} {\bf M} : {\bf K} + \rho^0|{\bf M}| |{\bf K}|]}, \quad x_3 = -\frac{1}{\rho^0 Re},\\
\\
y_{1,2} = \frac{- \frac{1}{{Re}} |{\bf M}| |{\bf K}| - x_{1,2}^2(\frac{1}{{Re}} + \rho^0 {\bf M}  : {\bf C})
- x_{1,2}( \frac{1}{{Re}} {\bf {\bf M}} : {\bf K} + \rho^0 |{\bf M}|d)  }{3x_{1,2}^2 \rho^0 + 2x_{1,2}\rho^0({\bf M} : {\bf K}) + \rho^0 |{\bf M}| |{\bf K}|},\\
\\
y_3 = -\frac{1}{\rho^0 |{\bf M}| |{\bf K}|} \big[ x_3^2(\rho^0 {\bf M} : {\bf K}) + x_3 (\rho^0 |{\bf M}|d+ \frac{1}{{Re}}  {\bf M} : {\bf K})
+ |{\bf M}| \frac{1}{{Re}} d + g_1 |{\bf K}| \big].
\eea\een
 The thermodynamic mode $\alpha_1$   is related to the mobility matrix and hessian matrix of the bulk free energy exclusively. The rest two eigenvalues  $\alpha_{2,3}$ are coupled with hydrodynamics.


\begin{table}
\caption{Sign of the eigenvalues when $|k| \ll 1$ in different regimes of {\bf C}. Negative sign indicates stability while positive sign indicates instability. }
\begin{center}
\begin{tabular}{|c|c|c|c|c|}
\hline
                          & $\alpha_0$  &  $\alpha_1$                      & $\alpha_2$  & $\alpha_3$  \\
                          \hline
 {\bf C} $>$ 0          & negative                 &    negative                                &    negative    &        negative                   \\
 \hline
{\bf C} $<$ 0          & negative                 &     positive                              &    positive         &         negative            \\
 \hline
$ {\bf C}$ is & negative                  &   ${\bf p}^T \cdot  {\bf C} \cdot {\bf p} $ has the same sign with  &If ${\bf p}^T \cdot  {\bf C} \cdot {\bf p} > 0$: negative;         & negative                    \\
  indefinite             &         & $| {\bf C}|$: negative;  Otherwise, positive.   &If ${\bf p}^T \cdot  {\bf C} \cdot {\bf p}< 0$: positive.       &             \\
\hline
\end{tabular}
\end{center}
\label{table:stability_category2}
\end{table}

 Obviously, $\alpha_0$ is negative so the viscous mode is stable. From the asymptotic expansions of $\alpha$ at $|k| \gg 1$, we observe that all three eigenvalues $\alpha_{1,2,3}$ are negative (\ref{eq:alpha_k_Inf_com1}) since ${\bf K} > 0, {\bf M} \geq 0$ and viscosity coefficients positive. This indicates that the model does not have any short-wave instability near its steady states, which is physically meaningful.

When $|k| \ll 1$, we notice  that the leading term in $\alpha_1$ is determined by the combination of mobility coefficient matrix ${\bf M}$ and hessian matrix ${\bf C}$ of the bulk free energy. We assume that ${\bf M} \geq 0$ and has at least one positive eigenvalue,  so $g_1 > 0$. We discuss the dependence of the leading order term of $\alpha_1$ on ${\bf C}$.
\begin{itemize}
\item  When ${\bf C}  > 0$,   the leading term $- \frac{g_1 |{\bf C}|}{{\bf p}^T \cdot  {\bf C} \cdot  {\bf p}} k^2 < 0$, then $\alpha_1 < 0$. So, this mode is stable.

\item When ${\bf C} < 0$, the leading term $ - \frac{g_1 |{\bf C}|}{{\bf p}^T \cdot {\bf C} \cdot  {\bf p}} k^2 > 0$, then $\alpha_1 > 0$. This instability is due to the spinodal decomposition in the coupled Cahn-Hilliard type equations of $\rho_1$ and $\rho_2$.

\item  When ${\bf C}$ is indefinite and ${\bf p}^T \cdot {\bf C} \cdot {\bf p} $ has the same sign with $|{\bf C}|$ , the property of $\alpha_{1}$ is the same as the case where ${\bf C} > 0$;  Otherwise, the property of $\alpha_{1}$  is the same as the case of ${\bf C} < 0$.
\end{itemize}
$\alpha_{2,3}$ represent the two coupled modes. Their signs depend on the model parameters. Since the leading term is determined by the properties of the hessian matrix ${\bf C}$, we discuss their dependence on ${\bf C}$ below.
\begin{itemize}
\item  When ${\bf C}  > 0$,   $\sqrt{(-\frac{1}{\rho^0} {\bf p}^T \cdot {\bf C} \cdot {\bf p})}$ is imaginary.   In this situation, the leading order growth rate in $\alpha_{2,3}$ is the quadratic term $(- \frac{1}{{Re}} \frac{1}{2 \rho^0} - \frac{1}{2{\bf p}^T \cdot {\bf C} \cdot {\bf p}}( M_{11}(\rho_1^0 h_{\rho_1 \rho_1} + \rho_2^0 h_{\rho_1 \rho_2})^2 + M_{22}(\rho_1^0 h_{\rho_1 \rho_2} + \rho_2^0 h_{\rho_2 \rho_2})^2) )k^2 \leq 0 $. So, the two modes are stable.

\item When ${\bf C} < 0$, the leading term is given by   $\pm \sqrt{(-\frac{1}{\rho^0} {\bf p}^T \cdot {\bf C} \cdot {\bf p})}k$, indicating there exists an  unstable mode. This verified the fact that the steady state at a concave  free energy surface is unstable.

\item  When ${\bf C}$ is indefinite and ${\bf p}^T \cdot {\bf C} \cdot {\bf p} > 0$, the property of $\alpha_{2,3}$ is the same as the case where ${\bf C} > 0$. Similarly, the property of $\alpha_{2,3}$  is the same as the case of ${\bf C} < 0$ when ${\bf p}^T \cdot {\bf C} \cdot {\bf p}< 0$.
\end{itemize}
The stability property of the model with respect to ${\bf C}$ in the long wave regime is summarized in Table (\ref{table:stability_category2}). For the intermediate wave regime, we have to compute the growth rate numerically, which can only be done for specific  free energy density functions.

\subsection{Compressible model with the local mass conservation law}

Notice that the compressible model with the local mass conservation law also admits the same constant solution \eqref{csol}.
We repeat the same normal mode analysis analogous to the previous model and obtain the   dispersion equation as follows:
\ben\bea{l}
(\alpha \rho^0 + \frac{1}{{Re}_s} k^2)  \{\alpha^3 \rho^0 + \alpha^2[\rho^0k^2M_{11}(  \tilde{h}_{\rho_1 \rho_1} +  \tilde{\kappa}_{\rho_1 \rho_1}k^2) + \frac{1}{{Re}} k^2]
+ \alpha[k^2M_{11}( \tilde{h}_{\rho_1\rho_1} +  \\
\\
\tilde{\kappa}_{\rho_1 \rho_1}k^2) \frac{1}{{Re}} k^2 +
 {\bf p}^T\cdot  {\bf C} \cdot {\bf p}k^2+{\bf p}^T \cdot  {\bf K} \cdot {\bf p}k^4))]
+ k^4 M_{11} (\rho^0)^2 (( \tilde{h}_{\rho_1 \rho_1} + \\
\\
k^2  \tilde{\kappa}_{\rho_1 \rho_1})( \tilde{h}_{\rho \rho} +  \tilde{\kappa}_{\rho \rho}k^2)
-( \tilde{h}_{\rho \rho_1} + k^2  \tilde{\kappa}_{\rho \rho_1})^2 )  \}
= 0.
\eea\een
Again, $\alpha_0 = - \frac{1}{\rho^0} \frac{1}{{Re}_s} k^2$ is a root of this algebraic equation. We present the rest   asymptotically.

When $|k| \ll 1$, we have
\ben\bea{l}
\alpha_{1} = - \frac{M_{11}  (\rho^0)^2 | {\bf C}|}{{\bf p}^T \cdot  {\bf C} \cdot {\bf p}}k^2 +  (-\frac{x_0^3 \rho^0 + x_0^2[\rho^0 M_{11}  \tilde{h}_{\rho_1 \rho_1} + \frac{1}{{Re}} ] + x_0[{\bf p}^T \cdot  {\bf K} \cdot {\bf p} + \tilde{h}_{\rho_1 \rho_1} M_{11} \frac{1}{{Re}} ]}{{\bf p}^T \cdot  {\bf C} \cdot  {\bf p}}\\\\
-\frac{M_{11} (\rho^0)^2  [ \tilde{h}_{\rho_1\rho_1}  \tilde{\kappa}_{\rho \rho} +  \tilde{h}_{\rho \rho}  \tilde{\kappa}_{\rho_1 \rho_1} - 2  \tilde{h}_{\rho \rho_1}  \tilde{\kappa}_{\rho \rho_1}]}{{\bf p}^T \cdot  {\bf C} \cdot {\bf p}}) k^4+ O(k^5),\\\\
\alpha_{2,3} = \pm  \sqrt{(-\frac{1}{\rho^0} {\bf p}^T\cdot   {\bf C} \cdot {\bf p})}k - (\frac{1}{{Re}} \frac{1}{2 \rho^0}+ \frac{ M_{11}}{2{\bf p}^T\cdot   {\bf C} \cdot {\bf p}}  (\rho_1^0  \tilde{h}_{\rho_1 \rho_1}+\rho^0  \tilde{h}_{\rho \rho_1})^2)k^2  +   O(k^3),
\eea\label{eq:alpha_k_small}
\een
where $x_0 =  - \frac{M_{11} (\rho^0)^2 | {\bf C}|}{{\bf p}^T \cdot  {\bf C} \cdot {\bf p}}$.
When $|k| \gg 1$,
\ben\bea{l}
\alpha_1 =- M_{11}  \tilde{\kappa}_{\rho_1 \rho_1} k^4 -  M_{11}  \tilde{h}_{\rho_1 \rho_1} k^2 + O(k),\\\\
\alpha_{2,3} =  \frac{- \frac{1}{{Re}} \pm \sqrt{( \frac{1}{{Re}} )^2 - 4   \tilde{\kappa}_{\rho_1 \rho_1}^{-1} (\rho^0)^3 |{\bf K}|}}{2  \rho^0} k^2 + \frac{-M_{11} (\rho^0)^2 [\tilde{h}_{\rho_1 \rho_1}  \tilde{\kappa}_{\rho \rho} +  \tilde{h}_{\rho \rho}   \tilde{\kappa}_{\rho_1 \rho_1}  - 2 \tilde{h}_{\rho \rho_1}  \tilde{\kappa}_{\rho \rho_1}]}{  2 x_{2,3}\rho^0 M_{11}  \tilde{\kappa}_{\rho_1 \rho_1} + M_{11}  \tilde{\kappa}_{\rho_1 \rho_1}\frac{1}{{Re}}  } \\\\
- \frac{x_{2,3}^3 \rho^0 + x_{2,3}^2 [\rho^0M_{11} \tilde{h}_{\rho_1 \rho_1} + \frac{1}{{Re}}] + x_{2,3}[M_{11} \tilde{h}_{\rho_1 \rho_1} \frac{1}{{Re}} + {\bf p}^T \cdot  {\bf K} \cdot {\bf p}]}{2 x_{2,3}\rho^0 M_{11}   \tilde{\kappa}_{\rho_1 \rho_1} + M_{11}  \tilde{\kappa}_{\rho_1 \rho_1} \frac{1}{{Re}}} + O(\frac{1}{k}),
\eea\label{eq:alpha_k_Inf}
\een
where $x_{2,3} = \frac{- \frac{1}{{Re}} \pm \sqrt{( \frac{1}{{Re}} )^2 - 4  \tilde{\kappa}_{\rho_1 \rho_1}^{-1} ( \rho^0)^3 |{\bf K}|}}{2  \rho^0}$,
\ben
 {\bf C} = \left(
\bea{cc}
 \tilde{h}_{\rho \rho}  & \tilde{h}_{\rho \rho_1}  \\
 \tilde{h}_{\rho \rho_1}  & \tilde{h}_{\rho_1 \rho_1}
\eea
\right)
\een
is the hessian matrix of the bulk free energy density function h with respect to  $\rho$ and $\rho_1$ and evaluated at the constant steady state, and
\ben
 {\bf K} = \left(
\bea{cc}
 \tilde{\kappa}_{\rho \rho}  &  \tilde{\kappa}_{\rho \rho_1}  \\
  \tilde{\kappa}_{\rho \rho_1}  &  \tilde{\kappa}_{\rho_1 \rho_1}
\eea
\right)
\een
is the coefficient matrix of the quadratic conformational entropy term in the free energy density function \eqref{free_energy_local}.

Like in the previous model, the first growth rate $\alpha_0$ is the viscous mode associated to the viscous stress exclusively; the second growth rate $\alpha_1$ is a thermodynamic mode, related to the transport equation of density $\rho_1$ and dictated by the mobility matrix and hessian matrix of the bulk free energy. The rest two growth rates  $\alpha_{2,3}$ are coupled modes.

 Obviously, $\alpha_0$ is negative so the viscous mode is stable. For the other three modes, we adopt the same strategy combining asymptotic analysis with numerical computations.
From asymptotic expansions  (\ref{eq:alpha_k_Inf}) of $\alpha$ at $|k| \gg 1$, we observe that all three eigenvalues $\alpha_{1,2,3}$ are negative, given that $ {\bf K}$ and the mobility coefficients  are both positive definite. This indicates that the model does not have any short-wave instability near its steady states. The properties of the three modes in the long wave regime are identical to the cases discussed in the previous section for the more general compressible model and summarized in Table \ref{table:stability_category2}.

For the intermediate wave regime, we have to compute the growth rate using a specific free energy density function numerically. We use the Peng-Robinson bulk free energy as an example here \cite{Peng_Robinson}, which is given by
\ben\bea{l}
\tilde{h}(\rho_1, \rho, T) =  \frac{r_m \rho_1 + (\rho - \rho_1) }{m_2} \varphi(T) -  \frac{r_m \rho_1 + (\rho - \rho_1)}{m_2} RT ln( \frac{m_2}{r_m \rho_1 + (\rho - \rho_1)} - b)
\\
-  \frac{r_m \rho_1 + (\rho - \rho_1)}{m_2}\frac{ a}{2\sqrt{2} b} ln[\frac{m_2 + (r_m \rho_1 + (\rho - \rho_1)) b(1+\sqrt{2})}{m_2 + (r_m \rho_1 + (\rho - \rho_1)) b(1-\sqrt{2})}] \\

\begin{figure}
\centering
\subfigure {\includegraphics[width=0.65\textwidth]{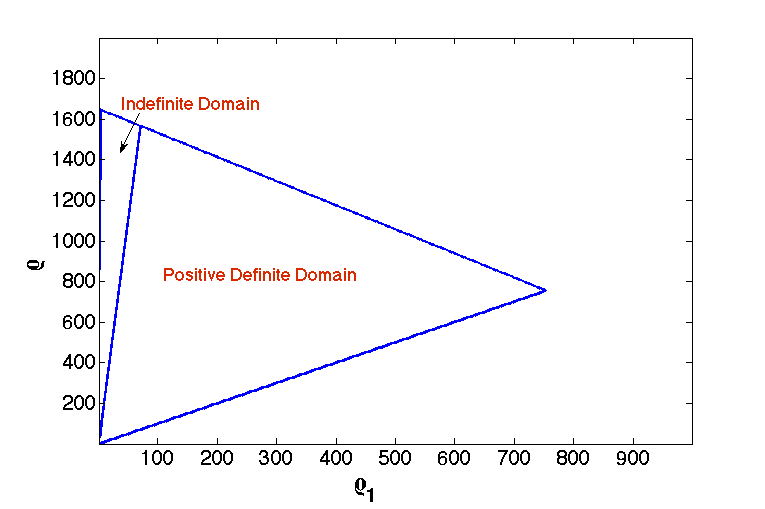}}\\
\caption{Domain of concavity of the Peng-Robinson free energy. }
\label{fig:Concavity}
\end{figure}

+  \frac{r_m \rho_1 + (\rho - \rho_1)}{m_2} RT[\frac{r_m \rho_1}{r_m \rho_1 + (\rho - \rho_1)}ln\frac{r_m \rho_1}{r_m \rho_1 + (\rho - \rho_1)}+ \frac{\rho_2}{r_m \rho_1 + (\rho - \rho_1)}ln\frac{(\rho - \rho_1)}{r_m \rho_1 + (\rho - \rho_1)}].
\eea\label{PR_free_energy2}
\een
This is obtained by replacing $\rho_1$, $\rho_2$ in the free energy density given in (\ref{PR_free_energy}) by $\rho_1, \rho_2=\rho-\rho_1$. This free energy density is either positive definite or indefinite in its entire physical domain. The positive definite domain and indefinite domain are shown in Figure \ref{fig:Concavity} in $(\rho_1, \rho)$ space. {\it Notice that in this example, when ${\bf C}$ is indefinite, we always have $|{\bf C}| < 0$, it is impossible to have two unstable modes $\alpha_1$ and $\alpha_2$ exist simultaneously according to table \eqref{table:stability_category2}.} We then search the parameter space to sample all the possible instabilities associated to the compressible model with this free energy.

\begin{figure}
\centering
\subfigure[ $\alpha_1$ when $|k| \ll 1$ ]{\includegraphics[width=0.325\textwidth]{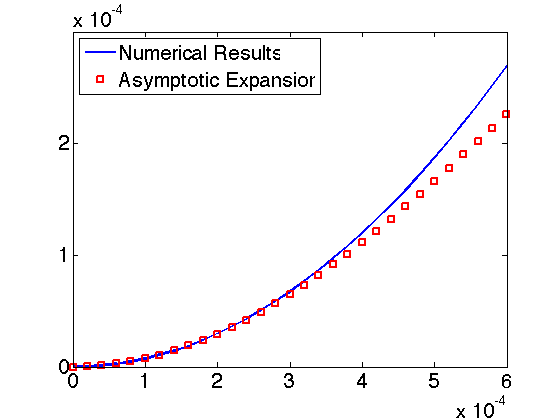}}
\subfigure[$\alpha_{1}$]{\includegraphics[width=0.325\textwidth]{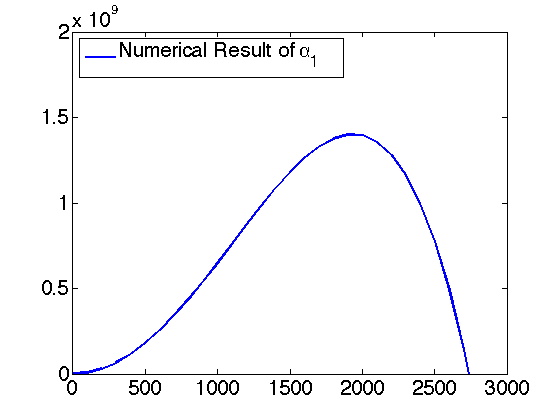}}
\subfigure[ $\alpha_1$ when $|k| \gg 1$ ]{\includegraphics[width=0.325\textwidth]{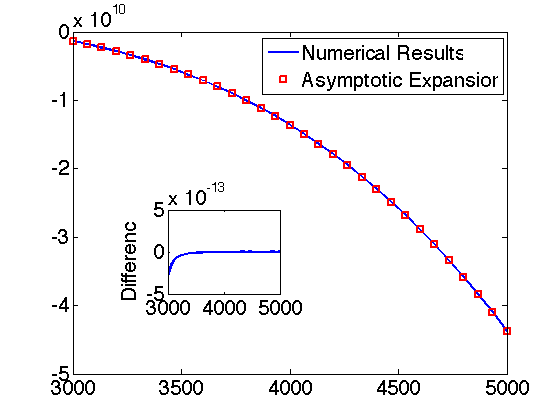}}\\
\subfigure[ $\alpha_{2,3}$ when $|k| \gg 1$ ]{\includegraphics[width=0.325\textwidth]{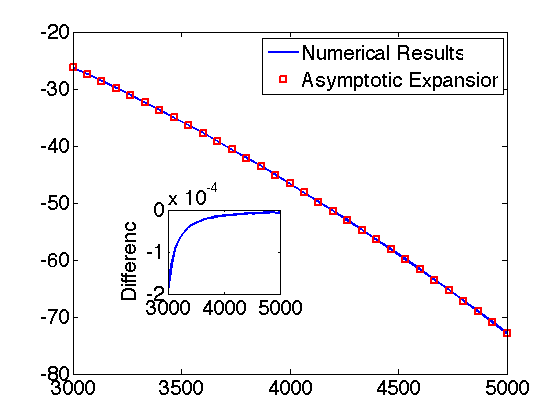}}
\subfigure[  $\alpha_{2,3}$ ]{\includegraphics[width=0.325\textwidth]{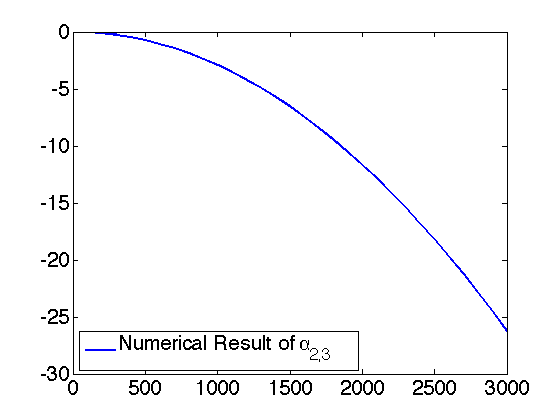}}
\subfigure[ $\alpha_{2,3}$ when $|k| \ll 1$ ]{\includegraphics[width=0.325\textwidth]{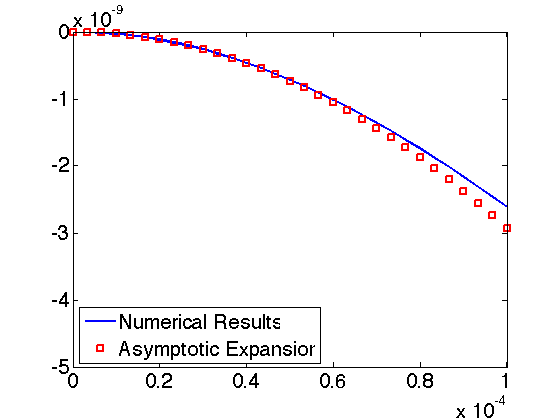}}

\caption{Numerical growth rates and the corresponding asymptotic ones as functions of the wave number when $\alpha_1>0$ and all others are negative in compressible model \eqref{eq:Nondim_com1} at constant state $(\rho^0 , \rho_1^0, {\bf v}) = (400, 2, 0 , 0)$ with the Peng-Robinson free energy. The vertical axis is the growth rate and the horizontal one is the wave number. (a). $\alpha_1$ in the long wave range. (b). $\alpha_1$ in the intermediate wave range. (c). $\alpha_1$
in the short wave range. (d). $\alpha_{2,3}$ in the short wave range. (e). $\alpha_{2,3}$ in the intermediate wave range. (f). $\alpha_{2,3}$ in the short wave range. The  parameter values used are: $M_{11} = 0.0001$, $Re_s =  1$, $Re_v = 3$, $ \tilde{\kappa}_{\rho \rho}= 0.000106$, $ \tilde{\kappa}_{\rho_1 \rho_1} = 0.0001$,  $ \tilde{\kappa}_{\rho \rho_1} = 0$. }
\label{fig:New_Asym_Num_case1}
\end{figure}

\begin{figure}
\centering
\subfigure[ $\alpha_{1}$ when $|k| \ll 1$ ]{\includegraphics[width=0.325\textwidth]{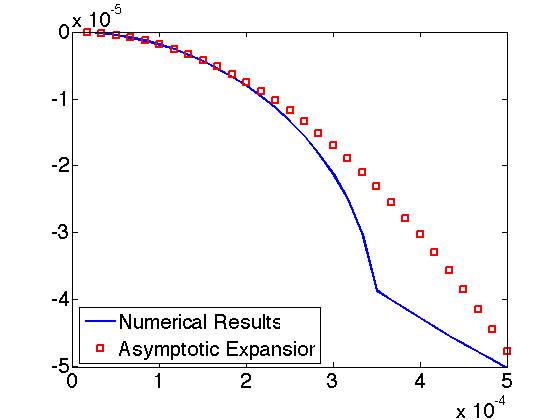}}
\subfigure[  $\alpha_{1}$ ]{\includegraphics[width=0.325\textwidth]{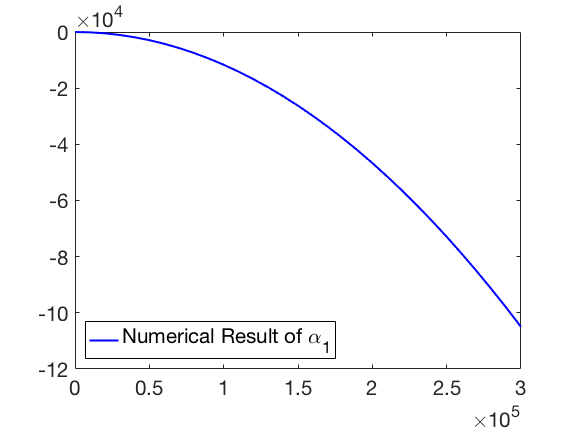}}
\subfigure[ $\alpha_{1}$ when $|k| \gg 1$ ]{\includegraphics[width=0.325\textwidth]{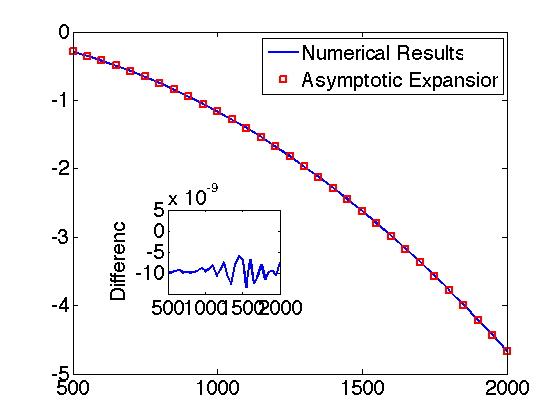}}\\
\subfigure[ $\alpha_2$ when $|k| \ll 1$ ]{\includegraphics[width=0.325\textwidth]{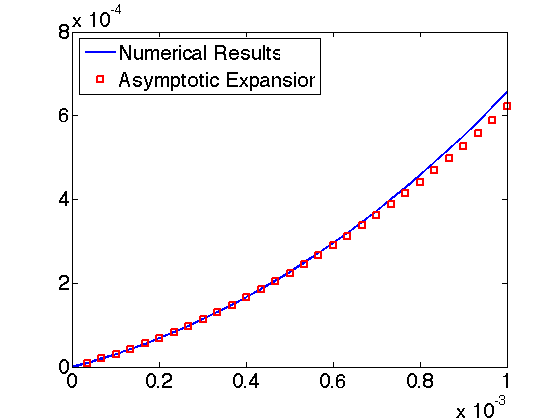}}
\subfigure[$\alpha_{2}$]{\includegraphics[width=0.325\textwidth]{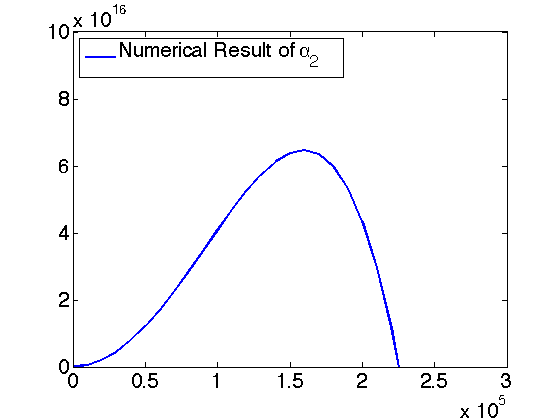}}
\subfigure[ $\alpha_2$ when $|k| \gg 1$ ]{\includegraphics[width=0.325\textwidth]{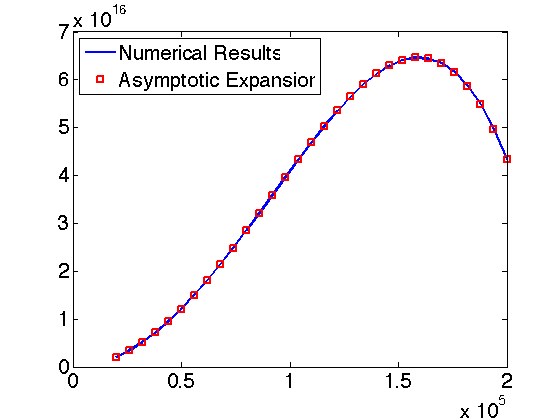}}\\
\subfigure[ $\alpha_{3}$ when $|k| \ll 1$ ]{\includegraphics[width=0.325\textwidth]{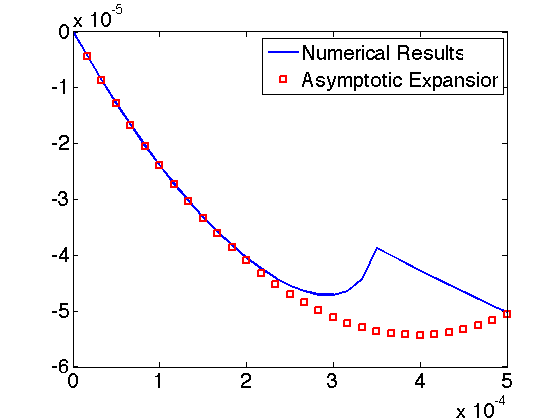}}
\subfigure[  $\alpha_{3}$ ]{\includegraphics[width=0.325\textwidth]{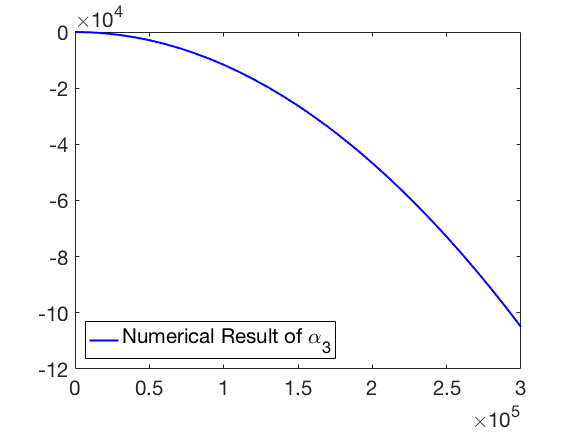}}
\subfigure[ $\alpha_{3}$ when $|k| \gg 1$ ]{\includegraphics[width=0.325\textwidth]{./fig5/case2_alpha13_k_inf}}
\caption{Numerical growth rates and the corresponding asymptotic ones when $\alpha_2>0$ while the others are negative compressible model \eqref{eq:Nondim_com1} at constant state $(\rho^0 , \rho_1^0, {\bf v}) = (1000, 0.025, 0 , 0)$ with the Peng-robinson free energy. The vertical axis is the growth rate and the horizontal one is the wave number.  (a). (d). (g). Growth rates in the short wave range. (b). (e). (h). Growth rates in the intermediate wave range. (c). (f). (j). Growth rates in the short wave range. The  parameter values used are:  $M_{11} = 0.0001$, $Re_s =  1$, $Re_v = 3$, $ \tilde{\kappa}_{\rho \rho}= 0.000106$, $ \tilde{\kappa}_{\rho_1 \rho_1} = 0.0001$,  $ \tilde{\kappa}_{\rho \rho_1} = 0$. }
\label{fig:New_Asym_Num_case2}
\end{figure}

\begin{figure}
\centering
\subfigure[ $\alpha_1$ when $|k| \ll 1$ ]{\includegraphics[width=0.325\textwidth]{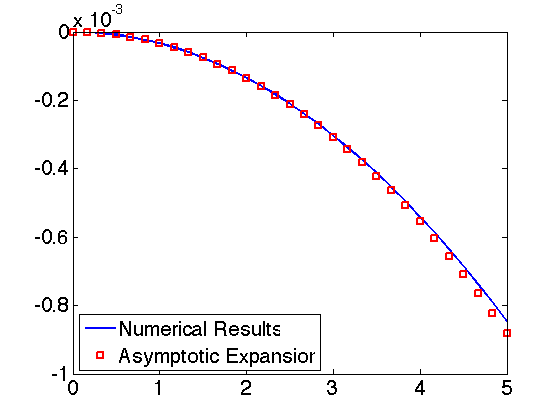}}
\subfigure[ $\alpha_{1}$  ]{\includegraphics[width=0.325\textwidth]{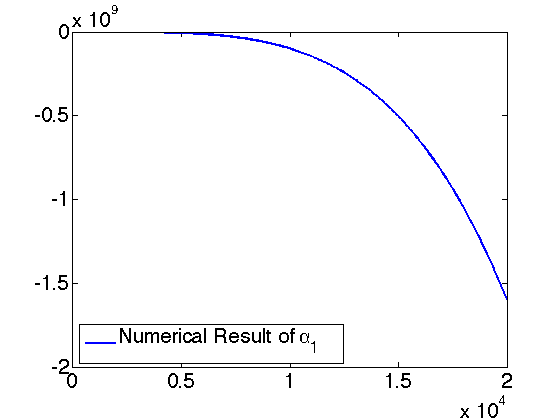}}
\subfigure[ $\alpha_{1}$ when $|k| \gg 1$ ]{\includegraphics[width=0.325\textwidth]{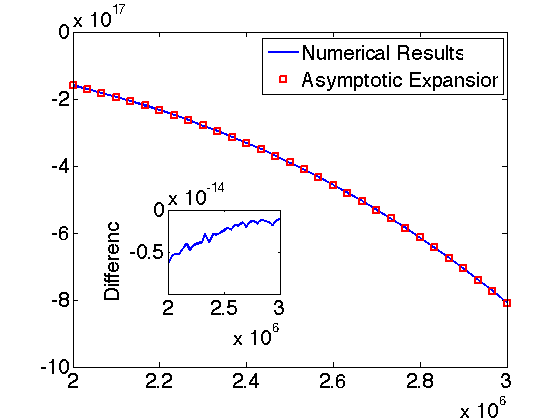}}\\
\subfigure[ $\alpha_{2,3}$ when $|k| \ll 1$ ]{\includegraphics[width=0.325\textwidth]{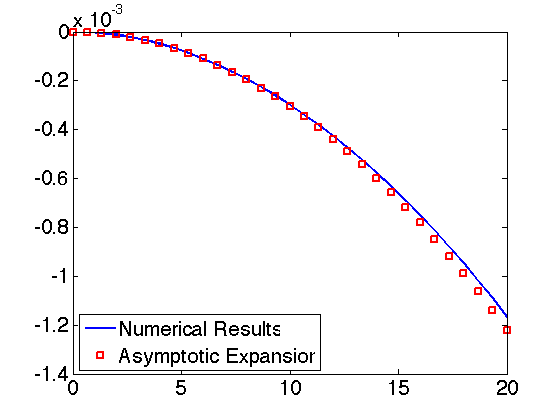}}
\subfigure[$\alpha_{2,3}$]{\includegraphics[width=0.325\textwidth]{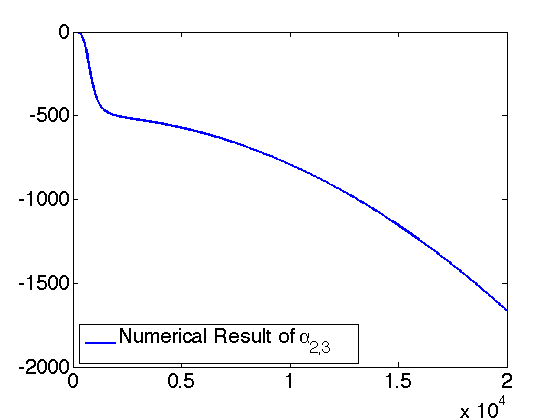}}
\subfigure[ $\alpha_{2,3}$ when $|k| \gg 1$ ]{\includegraphics[width=0.325\textwidth]{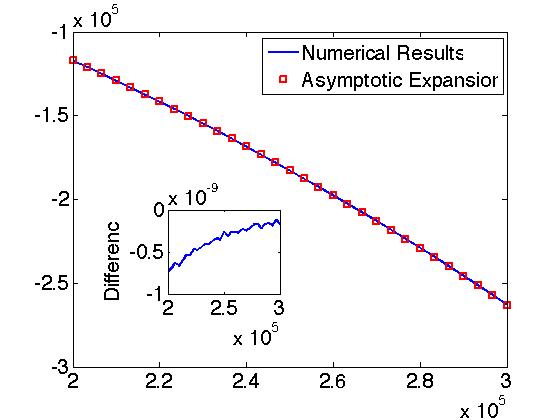}}\\
\caption{Numerical growth rates and the corresponding asymptotic ones without any unstable modes in compressible model \eqref{eq:Nondim_com1} at constant state $(\rho^0 , \rho_1^0, {\bf v}) = (400, 200, 0, 0)$ with the Peng-Robinson free energy. (a). and (d). Growth rates in the short wave range. (b). and (e).  Growth rates in the intermediate wave range. (c). and (f).  Growth rates in the short wave range. The parameter values used are: $M_{11} = 0.0001$, $Re_s =  10^{6}$, $Re_v = 3 \times 10^{6}$, $ \tilde{\kappa}_{\rho \rho}= 0.000106$, $ \tilde{\kappa}_{\rho_1 \rho_1} = 0.0001$,  $ \tilde{\kappa}_{\rho \rho_1} = 0$. }
\label{fig:New_Asym_Num_case3}
\end{figure}

As an example, we choose the steady state given by $(\rho^0 , \rho_1^0, {\bf v}_0) = (400, 2, 0, 0)$  to show the positive grow in $\alpha_1$. To show positive growth in the coupled mode $\alpha_2$, we choose $(\rho^0 , \rho_1^0, {\bf v}_0) = (1000, 0.025, 0, 0)$.  Figure \ref{fig:New_Asym_Num_case1} plots the three growth rates $\alpha_{1,2,3}$ with  $\alpha_1>0$ at the first constant solution. The corresponding eigenvector to $\alpha_1$ of the linearized system is $(0, 1, 0,0)$, indicating the unstable variable in the linear regime is $\rho_1$.  The three growth rates $\alpha_{1,2,3}$ with the coupled mode $\alpha_2>0$ at the second solution are plotted in Figure \ref{fig:New_Asym_Num_case2}.  The corresponding eigenvector to $\alpha_2$ is $(0, 1, 0 ,0)$ as well, indicating the instability is still associated with $\rho_1$. When $ {\bf C}>0$, the corresponding constant solution is stable. We choose constant solution $(\rho^0 , \rho_1^0, {\bf v}_0) = (400, 200, 0, 0) $ as an example. The three growth rates $\alpha_{1,2,3}$ of negative real parts are shown in Figure \ref{fig:New_Asym_Num_case3}. The numerical results show that the asymptotic analysis is  accurate in their respective wave number range of applicability.

From the linear analysis above, we conclude that linear dynamics of compressible model \eqref{eq:Nondim_com2} and  \eqref{eq:Nondim_com1} are  qualitatively the same. Next, we investigate the near equilibrium  dynamics of the quasi-incompressible model.

\subsection{Quasi-incompressible model}
The quasi-incompressible fluid flow model equations admit a constant solution:
\ben
{\bf v}={\bf 0},\quad \phi = \phi^0, \quad \Pi=\Pi_0,
\een
where $\phi^0,  \Pi_0$ are constants. We perturb the constant solution as follows:
\ben
{\bf v}= \epsilon e^{\alpha t + i {\bf k} \cdot {\bf x}} {\bf v}^{c}, \quad \phi=\phi^0 + \epsilon e^{\alpha t + i {\bf k} \cdot {\bf x}} {\phi}^{c},\quad \Pi=\Pi_0 + \epsilon e^{\alpha t + i {\bf k} \cdot {\bf x}} {\Pi}^{c},
\een
 where $\epsilon$ is a small perturbation, and ${\bf v}^c, \phi^c, \Pi^c$ are constants.

The dispersion equation is a factorable, third order polynomial in $\alpha$
\ben\bea{l}
(\frac{1}{{Re}_s} k^2 +\alpha \rho^0 )  (\alpha^2  (1- \frac{\hat{\rho_1}}{\hat{\rho_2}})^2\frac{1}{\hat{\rho_1}^2}M_{11} k^2 \rho^0
+ \alpha[ k^2 +
 \frac{1}{Re} (1-\frac{\hat{\rho_1}}{\hat{\rho_2}} )^2 \frac{1}{\hat{\rho_1}^2} M_{11}  k^4]\\
 \\
+     k^4M_{11} (  \hat{h}_{\phi \phi} + k^2 \hat{\kappa}_{\phi \phi} )  \frac{1}{\hat{\rho_1}^2}  [1 - (1-   \frac{\hat{\rho_1}}{\hat{\rho_2}}    )\phi_0 ]^2  ) = 0.
\eea\een
i.e.
\ben\bea{l}
(\frac{1}{{Re}_s} k^2 +\alpha \rho^0 )  [(1- \frac{\hat{\rho_1}}{\hat{\rho_2}})^2\frac{1}{\hat{\rho_1}^2}M_{11} k^2] (\alpha^2  \rho^0
+ \alpha[ [ (1- \frac{\hat{\rho_1}}{\hat{\rho_2}})^2\frac{1}{\hat{\rho_1}^2}M_{11} ]^{-1}  +
 \frac{1}{Re}  k^2]\\
 \\
+     k^2  (  \hat{h}_{\phi \phi} + k^2 \hat{\kappa}_{\phi \phi} )  [  \phi_0 -  \frac{ \hat{\rho_2}}{   \hat{\rho_2} - \hat{\rho_1} }   ]^2  ) = 0.
\eea\een
The growth rates are given explicitly by
\ben\bea{l}
\alpha_0 = - \frac{1}{{Re}_s} \frac{1}{\rho^0}k^2,\\
\\
\alpha_1 = \frac{- 2  k^2(\hat{h}_{\phi \phi} + k^2 \hat{\kappa}_{\phi \phi})Q^2}{[(\frac{1}{{Re}} k^2 + A )  + \sqrt{(\frac{1}{{Re}} k^2 + A )^2 - 4 \rho^0 k^2(\hat{h}_{\phi \phi} + k^2 \hat{\kappa}_{\phi \phi})Q^2}]},\\
\\
\alpha_2 = \frac{-(\frac{1}{{Re}} k^2 + A  )  - \sqrt{(\frac{1}{{Re}} k^2 + A )^2 - 4 \rho^0 k^2(  \hat{h}_{\phi \phi} + k^2 \hat{\kappa}_{\phi \phi}   )Q^2}}{2\rho^0},
\eea\een
where
\ben\bea{l}
Q = \phi^0 - \frac{\hat{\rho}_2}{\hat{\rho}_2 - \hat{\rho}_1}, \quad
\frac{1}{{Re}} = 2\frac{1}{{Re}_s} + \frac{1}{{Re}_v} > 0, \quad
A = [(1-\frac{\hat{\rho}_1}{\hat{\rho}_2})^2   \frac{1}{(\hat{\rho}_1)^2} M_{11}]^{-1} > 0 .
\eea\een
The stable hydrodynamic mode remains in $\alpha_0$. The thermodynamic modes are now given by $\alpha_{1,2}$. $Re(\alpha_1)$ can be positive only when  $\hat{h}_{\phi \phi} < 0$, in which $Re(\alpha_1)>0$  when $0 \leq k \leq \sqrt{-\frac{ \hat{h}_{\phi \phi}}{\kappa_{\phi \phi}}}$.  This instability is due to the spinodal decomposition in the coupled Cahn-Hilliard equation of $\phi$. Given that the viscosity and mobility coefficients are all positive, $Re(\alpha_2)<0$. So, the second coupled mode is a stable mode. In the long wave range ($|k|\ll 1$), $\alpha_1 \approx -\frac{M_{11} \hat{h}_{\phi \phi}}{\hat{\rho_1}^2\hat{\rho_2}^2}(\hat{\rho_1}+(\hat{\rho_1}-\hat{\rho_2})\phi^0)^2 k^2.$

When $\hat{\rho_1} = \hat{\rho_2}$, the model reduces to an incompressible model with the following two growth rates
\ben\bea{l}
\alpha_0 = - \frac{1}{{Re}_s} \frac{1}{\rho^0} k^2,\\
\\
\alpha_1 =  - \frac{1}{\hat{\rho_2}^2} M_{11} \hat{h}_{\phi \phi} k^2 - \frac{1}{\hat{\rho_1}^2} \hat{\kappa}_{\phi \phi} M_{11} k^4.
\eea \label{mode_incom}
\een
The thermodynamic mode decouples from the hydrodynamic mode completely in the linear regime. The possible instability only lies in the spinodal mode of the Cahn-Hilliard equation. In fact, $A, Q \to \infty $ in this limit. So, the growth rate associated with $\alpha_2$ in the quasi-incompressible model is lost.

\subsection{Summary of linear stability results}

 In compressible phase field models, there are four modes in the 1D perturbation analysis: $\alpha_0$ is the hydrodynamic mode dictated by the viscous stress, $\alpha_1$ is the thermodynamic mode dominated by the mobility and   the bulk free energy, the rest two modes $\alpha_{2,3}$ are coupled, which couples dynamics of phase behavior with hydrodynamics and may be unstable depending on the composition of the fluid mixture. When the Hessian matrix of the bulk free energy ${\bf C}>0$,   $\sqrt{(-\frac{1}{\rho^0} {\bf p}^T \cdot {\bf C} \cdot {\bf p})}$ is imaginary. So  $\pm \sqrt{(-\frac{1}{\rho^0} {\bf p}^T \cdot {\bf C} \cdot {\bf p})}k$ represents a wave that does not contribute to the amplitude change in growth rates of the linearized system. The scenario on stability of the steady state is tabulated in table \eqref{table:stability_category2}.

  When the quasi-incompressible constraint is added, i.e. $\rho_1 = \hat{\rho_1} \phi, \rho_2 = \hat{\rho_2} (1 - \phi)$. The positive definite matrix ${\bf C}$ reduces to a singular matrix
\ben
 {\bf C} =
h_{\phi \phi}  \left(\bea{cc}
 \frac{1}{\hat{\rho_1} \hat{\rho_1}} & - \frac{1}{\hat{\rho_1} \hat{\rho_2}}  \\
-  \frac{1}{\hat{\rho_1} \hat{\rho_2}}  &  \frac{1}{\hat{\rho_2} \hat{\rho_2}}
\eea
\right)
\een
Obviously, $|{\bf C}| = 0$ and ${\bf p}^T \cdot {\bf C} \cdot {\bf p} = (2\phi-1)^2$ for $\rho_1^0=\phi \hat{\rho_1}$, $\rho_2^0=(1-\phi) \hat{\rho_2}$. The growth rates reduce to two modes labeled as $\alpha_{1,2}$. They are not necessary related to the $\alpha_{1,2}$ in the compressible model.  Furthermore, when the quasi-incompressible mode reduces to the incompressible model, the coupled hydrodynamic modes vanishes, leading to one mode in $\alpha_1$.

The analysis shows that the more constraints we have on the composition of the fluid mixture, the less coupled the equations are in the linear regime.
 In 3D models, the total number of growth rates will increase as the number of equations increases. But, the number of unstable modes will not change. In addition, the 1D perturbation analysis in wave numbers applies to multi-dimensional case  as well. We will not omit the details for simplicity. 

\section{Conclusion}
We have presented a systematic way to derive hydrodynamic phase field models for  multi-component fluid mixtures of compressible fluids as well as incompressible fluids. The governing equations in the models are composed of the mass and momentum conservation law as well as the constitutive equations, which are derived using the generalized Onsager Principle to warrant an energy dissipation in time. By relaxing or enforcing local mass conservation law while keeping the total mass conserved, we obtain  two classes of  compressible models, one conserves the local mass while the other does not. Via a Lagrange multiplier approach, we reduce the compressible model with the local mass conservation law to a quasi-incompressible model when the constituent fluids are all incompressible. The quasi-incompressible model further reduces to the incompressible model.

We then study linear stability of all the models. The properties of linear stability are studied and differences of the models in the linear regime are identified: there exist three types of growth/decay rates among the models. The first type is dominated by the viscous property of the fluid, known as the viscous mode. The second type is the thermodynamic mode, which is dominated by the mobility and Hessian  of the bulk free energy density.  The third type is the coupled mode among the phase variables and hydrodynamic variables. When more constraints are enforced to reduce the models from the compressible, to the quasi-incompressible and then to the incompressible model, the number of coupled modes reduces accordingly, indicating that these constraints weaken the coupling of the equations in the model. This study not only develops a general framework for the derivation of compressible models and  their  reduction to quasi-incompressible models, but also identifies  differences between compressible and incompressible models in near equilibrium dynamics. It provides an easy to use theoretical tool for studying hydrodynamics of multiphasic fluids.

\section*{Acknowledgements}
Qi Wang's research is partially supported by NSF-DMS-1517347, DMS-1815921 and OIA-1655740, NSFC awards $\#$11571032, $\#$91630207 and NSAF-U1530401. Tiezheng Qian's research is partially supported by Hong Kong RGC Collaborative Research Fund No. C1018-17G.

\section{Appendix: Dispersion equations of the 2D hydrodynamic models}

We list the dispersion equations in determinant forms of all hydrodynamic models derived in this study in 2 space dimension in the appendix.

\subsection{Dispersion equation of the compressible model with the global mass conservation}
The dispersion equation of the linearized equation system of the compressible model with the global mass conservation is given by  a 4$\times$4 determinant as follows
\ben
\bea{l}
det
\left(
\bea{cccc}
\alpha +  A_{11}   &  A_{12}   &       i \rho_1^0  k   & 0            \\

 A_{21}  & \alpha + A_{22} &   i \rho_2^0  k     & 0        \\

ik (\rho_1^0 D_{11}  +  \rho_2^0  D_{12})   &  ik(\rho_2^0 D_{22} +  \rho_1^0  D_{12})      & \alpha \rho^0 +    \frac{1}{{Re}} k^2  & 0      \\
0 & 0  & 0 & \alpha \rho^0  + \frac{1}{{Re}_s} k^2        \\
\eea
\right) = 0
,
\eea\label{dispersion2}
\een
where $A_{11} = k^2 (M_{11}D_{11}   + M_{12}D_{12})$, $A_{12} = k^2 (M_{11}D_{12} +  M_{12}D_{22})$, $A_{21} =k^2 (M_{12}D_{11} + M_{22}D_{12} )  $, $A_{22} = k^2 (M_{12}D_{12}  +   M_{22}D_{22} ) $ and $D_{11} = h_{\rho_1 \rho_1} + k^2 \kappa_{\rho_1 \rho_1}$, $D_{22} = h_{\rho_2 \rho_2} + k^2 \kappa_{\rho_2 \rho_2}$, $D_{12} = h_{\rho_1 \rho_2} + k^2 \kappa_{\rho_1 \rho_2}$, $\frac{1}{{Re}} = 2\frac{1}{{Re}_s} + \frac{1}{{Re}_v}$. The growth/decay rate in the hydrodynamic mode associated to the viscous stress is given explicitly by  $\alpha=- \frac{1}{{Re}_s} \frac{1}{\rho^0} k^2$, which decouples from the rest of the modes. This decoupling is inherited by  all its limiting models given below.

\subsection{Dispersion equation of the compressible model with local mass conservation}
The dispersion equation of the linearized equation system of this model is given by  a 4$\times$4 determinant as follows
\ben
\bea{l}
det\left(
\bea{cccc}
\alpha     &   0   &   i \rho^0  k & 0      \\

( k^2 M_{11})D_{12}    & \alpha +  (k^2 M_{11})D_{22} &     i \rho_1^0  k    & 0       \\

ik(\rho_1^0  D_{12}  +  \rho^0  D_{11} )   &  ik(\rho_1^0 D_{22}+  \rho^0 D_{12})     & \alpha \rho^0 + \frac{1}{{Re}} k^2  & 0      \\
0 & 0 & 0 & \alpha \rho^0  +    \frac{1}{{Re}_s} k^2   \\
\eea
\right)=0,
,
\eea\label{dispersion1}
\een
where $D_{11} = \tilde{h}_{\rho \rho} + k^2  \tilde{\kappa}_{\rho \rho}$, $D_{22} = \tilde{h}_{\rho_1 \rho_1} + k^2   \tilde{\kappa}_{\rho_1 \rho_1}$, $D_{12} =   \tilde{h}_{\rho \rho_1} + k^2  \tilde{\kappa}_{\rho \rho_1}$, and $\frac{1}{{Re}} = 2\frac{1}{{Re}_s} +  \frac{1}{{Re}_v}$.

\subsection{Dispersion equation of the quasi-incompressible model}

The resulting dispersion equation of the linearized system of this model is given by a $4 \times 4$ determinant as follows
\ben
\bea{l}
det\left(
\bea{cccc}
0   &  - \alpha ( 1-\frac{\hat{\rho_1}}{\hat{\rho_2}} )    & ik - i k \phi^0 (1-\frac{\hat{\rho_1}}{\hat{\rho_2}}) & 0 \\
\frac{1}{\hat{\rho_1}^2 }M_{11} k^2 (1-\frac{\hat{\rho_1}}{\hat{\rho_2}})        & \alpha  + \frac{1}{\hat{\rho_1}^2}M_{11} k^2D_{\phi}   &  ik \phi^0 &    0      \\
ik & ik \phi^0D_{\phi} & \alpha \rho^0 + \frac{1}{{Re}} k^2  & 0 \\
  0  & 0   &  0 &\alpha \rho^0 + \frac{1}{{Re}_s} k^2
\eea
\right)=0,
\eea
\een
where $D_{\phi} = \hat{h}_{\phi \phi} + \hat{\kappa}_{\phi \phi}k^2$, $ \hat{h}_{\phi \phi} = \frac{\partial^2 h}{\partial \phi^2}$ is the second order derivative of the bulk free energy density function h with respect to volume fraction $\phi$ at the constant solution, and $\hat{\kappa}_{\phi \phi}$ is the coefficient of the conformational entropy.
If we multiply $( 1-\frac{\hat{\rho_1}}{\hat{\rho_2}} )$ by the second row and add it to the first row of the dispersion relation matrix, we obtain
\ben
\bea{l}
det\left(
\bea{cccc}
\frac{1}{\hat{\rho_1}^2}M_{11} (1- \frac{\hat{\rho_1}}{\hat{\rho_2}})^2 k^2    & (1- \frac{\hat{\rho_1}}{\hat{\rho_2}})\frac{1}{\hat{\rho_1}^2}M_{11} k^2D_{\phi}     & ik & 0 \\
\frac{1}{\hat{\rho_1}^2 }M_{11}  (1-\frac{\hat{\rho_1}}{\hat{\rho_2}})  k^2      & \alpha  + \frac{1}{\hat{\rho_1}^2}M_{11} k^2D_{\phi}   &  ik \phi^0 &    0      \\
ik & ik \phi^0D_{\phi} & \alpha \rho^0 + \frac{1}{{Re}} k^2  & 0 \\
  0  & 0   &  0 &\alpha \rho^0 + \frac{1}{{Re}_s} k^2
\eea
\right)=0,
\eea\label{disp-qi}
\een

\subsection{Dispersion equation of the incompressible model}
The dispersion equation of the linearized system of the incompressible model is given by
\ben
\bea{l}
det\left(
\bea{cccc}
0   & 0     & ik & 0 \\
0       & \alpha  + \frac{1}{\hat{\rho_1}^2}M_{11} k^2D_{\phi}   &  ik \phi^0 &    0      \\
ik & ik \phi^0D_{\phi} & \alpha \rho^0 + \frac{1}{{Re}} k^2  & 0 \\
  0  & 0   &  0 &\alpha \rho^0 + \frac{1}{{Re}_s} k^2
\eea
\right)=0.
\eea
\een
This can be obtain from that in the quasi-incompressible model by equating $\hat{\rho_1}=\hat{\rho_2}$ in \eqref{disp-qi}. 
\bibliographystyle{plain}
\bibliography{xp_thesis}{}

\end{document}